\documentclass[11pt,a4paper]{article}
\usepackage{latexsym}
\usepackage{graphicx}
\usepackage{epsfig}
\usepackage{amssymb}
\usepackage{amsmath}
\usepackage{epstopdf}
\usepackage{eufrak}
\usepackage{epstopdf}
\usepackage{pb-diagram}
\usepackage{url}

\usepackage[hmargin=2.5cm, vmargin=2.5cm]{geometry} % Margins

\newcommand\Beq{\begin{eqnarray}} 
\newcommand\Eeq{\end{eqnarray}}

\newcommand{\eq}[1]{equation~(\ref{#1})}
\newcommand{\eqs}[2]{equations~(\ref{#1})~\&~(\ref{#2})}
\newcommand{\eqss}[2]{equations~(\ref{#1})--(\ref{#2})}
\newcommand{\Eq}[1]{Equation~(\ref{#1})}
\newcommand{\Eqs}[2]{Equations~(\ref{#1})~\&~(\ref{#2})}
\newcommand{\Eqss}[2]{Equations~(\ref{#1})--(\ref{#2})}

\newcommand{\ie}{\textit{i.e.}, }
\newcommand{\eg}{\textit{e.g.}, }

\newcommand{\Reg}[1]{\mathrm{Reg}{(\ell#1)}}

\newcommand{\dd}[1]{\,\mathrm{d}{#1}}

\newcommand{\D}[1]{\nabla_{\! #1}}
\newcommand{\e}[1]{e_{#1}}

\newcommand{\pd}[1]{ \partial_{#1} }

\newcommand{\z}{z}

\begin{document}

\renewcommand*{\thefootnote}{\fnsymbol{footnote}}

\centerline{\textbf{\Large{Tensor calculus in spherical coordinates using Jacobi polynomials.}}}
\vspace{0.2cm}
\centerline{\textbf{\Large{Part-I: Mathematical analysis and derivations}}} 
\medskip
\centerline{\large{Geoff Vasil$^{*,1}$, Daniel Lecoanet$^{*,2,3}$, Keaton Burns$^{4}$, Jeff Oishi$^{5}$, Ben Brown$^{6}$}}
\medskip
\centerline{\textsl{\small{*Corresponding authors; email: geoffrey.vasil@sydney.edu.au; lecoanet@princeton.edu}}}
\smallskip
\centerline{\textsl{\small{$^1$University of Sydney School of Mathematics and Statistics, Sydney, NSW 2006, Australia}}}
\smallskip
\centerline{\textsl{\small{$^2$Princeton Center for Theoretical Science, Princeton, NJ 08544, USA}}}
\centerline{\textsl{\small{$^3$Princeton University Department of Astrophysical Sciences, Princeton, NJ 08544, USA}}}
\smallskip
\centerline{\textsl{\small{$^4$Massachusetts Institute of Technology Department of Physics, Cambridge, MA 02139, USA}}}
\smallskip
\centerline{\textsl{\small{$^5$Bates College Department of Physics and Astronomy, Lewiston, ME 04240, USA}}}
\smallskip
\centerline{\textsl{\small{$^6$University of Colorado Laboratory for Atmospheric and Space Physics and Department of}}}
\centerline{\textsl{\small{ Astrophysical and Planetary Sciences, Boulder, CO 80309, USA}}}

\renewcommand*{\thefootnote}{\arabic{footnote}}

\bigskip
\centerline{\large \today}
\bigskip

\bigskip
\textbf{Abstract}

This paper presents a method for the accurate and efficient computations on scalar, vector and tensor fields in three-dimensional spherical polar coordinates. The methods uses spin-weighted spherical harmonics in the angular directions and rescaled Jacobi polynomials in the radial direction. For the 2-sphere, spin-weighted harmonics allow for automating calculations in a fashion as similar to Fourier series as possible. Derivative operators act as wavenumber multiplication on a set of spectral coefficients. After transforming the angular directions, a set of orthogonal tensor rotations put the radially dependent spectral coefficients into individual spaces each obeying a particular regularity condition at the origin. These regularity spaces have remarkably simple properties under standard vector-calculus operations, such as \textit{grad} and $\textit{div}$. We use a hierarchy  of rescaled Jacobi polynomials for a basis on these regularity spaces. It is possible to select the Jacobi-polynomial parameters such that all relevant operators act in a minimally banded way. Altogether, the geometric structure allows for the accurate and efficient solution of general partial differential equations in the unit ball. 

\bigskip

\textbf{Keywords:} Spherical Geometry; Spectral Methods; Spin-weighted Spherical Harmonics; Jacobi Polynomials; Sparse Operators; 

\bigskip
\bigskip
\bigskip
\bigskip
\bigskip
\bigskip
\bigskip
\bigskip
\bigskip
\bigskip
\bigskip

\pagebreak

\section{Introduction}

This paper outlines a set of algorithms for arbitrary vector and tensor computations in spherical polar coordinates in the three-dimensional unit ball. A previous paper (Vasil \textit{et al.} 2016 \cite{V16}, from now on V16) presented a similar calculus for vector and tensor operations in cylindrical polar coordinates on the unit disk. The two-dimensional disk is one of the simplest examples of a well-behaved domain with a natural coordinate singularity. The unit ball has three separate coordinate singularities: disk-like singularities at the north and south pole, and a third at the centre of the domain.

Spherical geometry is important for a large number of two- and three-dimensional applications; see \eg \cite{Beyer_etal,Slevinsky_etal,townsend_etal_2016,Rubinstein,Libsharp,SHTns,Marti-benchmark,Hollerbach+2013,Kostelec}. Astrophysical and planetary applications provide many obvious examples with stars and planets. In this paper, we first discuss the mathematical structure for accurate tensor calculations in spherical coordinates. A companion paper (Part-II) discusses some applied test problems that demonstrate effectiveness in relevant geophysical and astrophysical scenarios.  These kinds of multi-scale problems often require high spatial resolution, but accurately representing functions in the vicinity of coordinate singularities can become numerically unstable in the presence of small-scale features. The most reliable way around this difficulty is to design a numerical method that explicitly conforms to the correct type of singular behaviour. 

A considerable amount of past work is devoted to overcoming problems with coordinate singularities in spheres and disks. Boyd \& Yu (2011) \cite{boyd2011comparing} provide an excellent review of the topic and compare several different methods. We also summarise the history and a number of the issues in V16. A Fourier series in polar coordinates demonstrates the issue most clearly. 
For example, 
\Beq
f(R,\phi) \ = \ \sum_{m=-\infty}^{\infty} f_{m}(R) e^{i m \phi}, \quad f_{m}(R) \ = \ \frac{1}{2\pi}\int_{0}^{2\pi} f(R,\phi) e^{-i m \phi} \dd{\phi},
\Eeq
where $R$ is the cylindrical radius and $\phi$ the circular angle.  If we assume $f$ is infinitely differentiable in the vicinity of the origin, then  
\Beq
\label{analytic-condition} f_{m}(R) \ \sim \ R^{|m|} F(R^{2} ) \quad \mathrm{as} \quad R \ \to \ 0,
\Eeq
where $F$ is a well-behaved function of $R^{2}$ \cite{orszag_1974}.  In V16 we called the radial and angular coordinates $(r,\theta)$.  To avoid confusion with 3D polar coordinates, we relabel the 2D coordinates $(r,\theta)  \to  (R , \phi)$. From now on, $0\le r\le1$ represents a 3D radius, $0\le \theta \le \pi$ represents a polar angle (co-latitude), and $0 \le \phi \le 2 \pi$ represents a circular angle (longitude). The circular angle is a common feature between the disk and the sphere. 

Representing a highly multi-scale function requires components with $|m| \gg 1$. Expanding $f_{m}(R)$ in a series of functions without explicitly considering the behaviour at the origin requires a large degree of cancelation to produce an accurate result.  Several authors discuss the numerical challenges of resolving \eq{analytic-condition} when $m \gg 1$ \cite{boyd_book,
boyd2011comparing,
Kostelec,
livermore_jones_worland_2007,
matsushima1995spectral,
sakai2009application,
zernike}.
The two-dimensional disk is a prototype for this kind of general behaviour. The polar regions in two-dimensional spherical coordinates are locally identical to the disk, with colatitude replacing the cylindrical radius. Fortunately, the 2-sphere comes equipped with spherical harmonic functions. In this case 
\Beq
\label{Y behaviour}
Y_{\ell,m}(\theta, \phi) \ \sim \ \sin^{|m|}(\theta) P(\cos (\theta)) e^{i m \phi},
\Eeq
where $P(\cos \theta)$ is a particular polynomial of $\cos \theta$. As $\theta \to 0$, $\sin( \theta) \sim \theta$ and $\cos (\theta) \sim 1 - \theta^{2}/2$; similar behaviours exist as $\theta \to \pi$. \Eq{Y behaviour} behaves exactly as \eq{analytic-condition} dictates it should, and spherical harmonic functions naturally conform to the singular behaviour of the poles at $\theta=0,\pi$. Therefore spherical harmonic functions naturally conform to the singular-coordinate behaviour at both the north and south poles.

With spherical harmonic functions we can represent analytic functions in the angular direction without fear from polar conditions,  
\Beq
\label{spherical scalar series}
f(r,\theta, \phi) &=& \sum_{m=-\infty}^{\infty}\sum_{\ell=|m|}^{\infty} f_{\ell,m}(r) Y_{\ell,m}(\theta, \phi), \\  
f_{\ell,m}(r) &=&  \int_{0}^{2\pi} \int_{0}^{\pi} f(r,\theta, \phi) {Y}_{\ell,m}(\theta,-\phi) \sin(\theta) \dd{\theta} \dd{\phi},
 \Eeq
where ${Y}_{\ell,m}(\theta,-\phi)$ gives the complex conjugate of ${Y}_{\ell,m}(\theta,\phi)$.

An analogous condition exists at the centre of the ball, 
\Beq
f_{\ell,m}(r) \ \sim \ r^{\ell} F(r^{2} ) \quad \mathrm{as} \quad r \ \to \ 0, \label{scalar reg.}
\Eeq
where $F$ is a well-behaved function of $r^{2}$. Throughout this paper, $\ell$ represents the spherical harmonic degree.  For components with $\ell\gg 1$, resolving behaviour at the centre of the sphere presents the same kind of numerical challenges as the centre of the disk. 

Scalar functions alone, like $f$ in equation~(\ref{scalar reg.}), are not sufficient for many important applications. In the vicinity of coordinate singularities, vector and tensor fields behave slightly differently than scalar functions. In the case of vectors, we know that the leading-order behaviour near the origin is at least one-degree lower than for scalar functions. For example, with the radial component of the gradient, 
\Beq
\hat{r} \cdot \nabla f \ \sim \  \frac{\pd{}}{\pd{} r} f_{\ell,m}(r) \ \sim \ r^{\ell-1}, \quad \mathrm{as} \quad r \ \to \ 0. \label{grad-r scaling}
\Eeq 
The angular gradient terms also multiply by $1/r$ terms that produce the same effect as \eq{grad-r scaling}, with terms scaling as $r^{\ell-1}$ as compared to the $r^{\ell}$ scaling of scalar functions~(\ref{scalar reg.}).

\textsl{\underline{Summary of results for the unit disk}: } In V16, we devised an efficient way to define and operate on functions, vectors and arbitrary tensors on the unit disk without problems resulting from the coordinate singularity at the origin.  Understanding the principles involved with the disk provide an excellent framework for understanding the more complicated calculus in three dimensional spheres. The main results regarding the disk are:
\begin{itemize}

\item 
The Fourier transform of vector-field components in the coordinate basis, $v^{R}_{m}(R)$, $v^{\phi}_{m}(R)$, naturally decompose into the sum of two alternative functions of radius, $v^{-}_{m}(R)$,  $v^{+}_{m}(R)$. One component scales near the origin as $v^{-}_{m}(R) \, \sim\, R^{|m|-1}$, the other scales as $v^{+}_{m}(R) \, \sim \, R^{|m|+1}$. Specifically,
$$v^{\pm}_{m}(R) \ \propto \  v_{m}^{R}(R) \pm i \,v_{m}^{\phi}(R) \ \sim \ R^{|m|\pm1} V^{\pm}(R^2) , \quad \mathrm{as} \quad R \ \to \ 0,$$
where $V^{\pm}$ are independent well-behaved functions of $R^{2}$. Higher-order tensor fields behave similarly.

\item
Given the regularity structure of vector and tensor fields, we used modified covariant derivative operators of the form 
$$\nabla_{\pm} \ \propto \  \nabla_{\!R} \pm i \nabla_{\!\phi}.$$
These derivative operators raise and/or lower the regularity class of tensor components. For the gradient of a scalar, $v^{\pm}(R) \propto  \nabla_{\pm} f(R)$. 

\item 
Given the above geometrical properties, V16 represented general tensor field components as a sum of functions of the form $$R^{|m|+\mu}\, P(2R^{2}-1)\, e^{i m \phi},$$ where $\mu$ is an integer that depends on the tensor rank.  $P(\z)$ is a polynomial (in the variable  $\z \equiv 2R^{2}-1$) carefully chosen to conform with the leading-order behaviour.  

\item 
It happens that Jacobi polynomials with parameters depending on $m$ and the tensor rank match the geometric regularity conditions in a natural way. The raising and lowering gradient operators increment and decrement the parameters of the Jacobi polynomials to match the geometrical regularity conditions of tensor components. The difference between past work and V16 was that we allowed the polynomial parameters to change under differentiation. 

\item 
A past solution for treating vectors was to use rescaled components and derivative operators such as $R \, v^{R}$, $R\, v^{\phi}$, $R\, \nabla_{\!R}$, $R\, \nabla_{\!\phi}$ \cite{sakai2009application}. 
While variables and operators of this type allow for sparse computations, they also have the effect of obscuring the underlying geometrical structure, which is somewhat of a nuisance. The main point of V16 was to exploit the underling geometrical structure in polar coordinates. This allowed for straightforward construction of complicated vector equations. The resulting matrix systems are even sparser than previous approaches, and better conditioned than what results from an approach that does not take the geometry into consideration.
  
\end{itemize}

In the current paper, we are interested in arbitrary vector and tensor calculations in spherical-polar coordinates. The same general ideas apply to the unit ball as the unit disk. The calculus of the unit 2-sphere is very similar to the unit disk. In that case, a complex-valued basis transformation produces vector components with a coherent \textit{spin-weight}. These components naturally break into a basis of \textit{spin-weighted spherical harmonics}. Several authors have used this approach on a number of different applications. Gelfand, \& Shapiro (1956) \cite{G_S_56} originally discovered spin-weighted spherical harmonics in their studies of the Lorenz group.  Newman \& Penrose (1966) \cite{N_P_66} applied the methods to classical problems in theoretical physics such as electro-magnetism and gravitational radiation.  Phinney \& Buridge (1973) \cite{P&B73} proposed (what they called) \textit{generalised spherical harmonics} as efficient analytical and computational tools in whole-earth geophysics. Our work relies on spin-weighted spherical harmonics as a necessary part of computations in three-dimensional spherical coordinates. 

In \S2 we review the most important features of tensor calculus on the unit 2-sphere, focusing on the essential building blocks for computations in the full 3-dimensional ball. Most of the analysis follows along the lines of \cite{P&B73}. We provide a few new results needed to make arbitrary numerical calculations; particularly we discuss operators for multiplying by functions of $\cos(\theta)$ and/or $\sin(\theta)$.  In addition, the calculus of spin-weighted spherical harmonic functions are not widely known in the fluid mechanics literature. It is worth reviewing for this reason alone. 

In \S3 we proceed where Phinney \& Buridge (1973) \cite{P&B73} left off. In Appendix C of that work, the authors provided explicit formulae for the covariant differentiation of scalars and vectors in 3D. The formulae show the radial dependence of various spin-weighted components. The most obvious thing to say about the results are that they are complicated. The same level of organisation does not seem to appear with the radial dependence as it does with the angular components. We show that there is a way to organise calculations, but that the operators must occur in spectral space. That is, there exist transformations that render vector calculus operations into sparse decoupled operations. But now the transformations depend on the spherical harmonic degree, $\ell$. The crux in the three-dimensional calculations is that the individual (radially dependent) spherical harmonic coefficients decompose into a superposition of functions with different regularity near the origin. 

In the case of scalar functions, we use the shorthand 
\Beq
f_{\ell}(r) \ \in \ \Reg{}. \label{Reg ell}
\Eeq
This assumes that $\ell$ is a non-negative integer. 
We wish for \eq{Reg ell} to be understood to mean exactly the same thing as \eq{scalar reg.}. That is, we define the vector space of functions 
\Beq
\Reg{} \ \equiv \ \left\{\, f\, : \, [0,1] \, \to \, \mathbb{C} \quad \mathrm{s.t.} \quad f(r) \ \sim \ r^{\ell} \, F(r^{2}) \quad \mathrm{as} \quad r \ \to \ 0 \, \right\},
\Eeq
where $F(r^{2})$ represents any even function of $r$ that is analytic in a neighbourhood of $r=0$. We can put any additional structure on the function $F(r^2)$ that we need to make further computations make sense.  We call the spaces $\Reg{}$ \textit{regularity} classes of degree $\ell$. It is important to note that 
\Beq
\Reg{+2} \ \subset \Reg{}. \label{R +2 in R}
\Eeq 
In the case of vectors, the spherical harmonic coefficients now have three ($\ell$-dependent) components. \Eq{grad-r scaling} implies that we need to consider $\Reg{-1}$ at the very least. One of the main findings of the current paper is that vectors break into the direct sum of regularity classes, 
\Beq
v_{\ell}(r) \ \in \ \Reg{-1} \oplus \Reg{} \oplus \Reg{+1}.
\Eeq
Continuing on from vectors, rank-2 tensor fields have nine ($\ell$-dependent) components. These components reside in a direct sum of $\mathrm{Reg}(j)$ where $\ell-2 \le j \le \ell+2$. This only represents five seperate spaces, so some of these spaces appear more than once to represent nine-component tensors. How exactly this happens follows directly from our derivations in \S3. 

For the case of the two-dimensional disk, we found in V16 that 
\Beq
v_{m}(R) \ \in \ \mathrm{Reg}(|m|-1) \oplus \mathrm{Reg}(|m|+1); \label{disk reg}
\Eeq
although this is not how it was presented at the time. However, writing vectors in this form highlights an important difference between the 2-disk and 3-ball. For the disk,
\Beq\label{disk embedding}
R \, v_{m}(R) \ \in \ \mathrm{Reg}(|m|) \oplus \mathrm{Reg}(|m|+2) \ \subset \ \mathrm{Reg}(|m|) \oplus  \mathrm{Reg}(|m|) ;
\Eeq
That is, it is possible to rescale vectors in two dimensions to behave the same as scalars at the origin. This allows a wider range of computational methods for the disk, \eg \cite{sakai2009application}. The main result of V16 is that although \eq{disk embedding} holds, it is better to exploit the structure of \eq{disk reg}. But the sphere does not allow the same kind of rescaling in radius. We need to build everything from scratch. 

In \S3, we show why a decomposition into regularity spaces emerges naturally from vector and tensor calculus. We find that the classic operators of vector calculus decompose into operators that map between the regularity spaces. But at the same time, we also find that the transformations between the components of vectors and tensors in coordinate bases and the regularity spaces is not extremely obvious. 

There has been some related work on tensor spherical harmonics in the past. The most important predecessor to our work that we know of is James (1976) \cite{James76}; we refer to this work as J76 from now on. This remarkable paper works out a full theory of tensor spherical harmonics for general-order tensors. J76 provides recursion formulae for constructing higher-order tensors out of lower-rank tensors. These formulae are based on the Wigner 3-$j$ coefficients. The paper also provides formulae for computing vector calculus operations (\textit{div, grad, curl, Laplacian}). However, J76 made most of the derivations using spherical coordinates for the angle variables and \textit{Cartesian} basis vectors for the tensorial behaviour. There is only a small discussion of the results using polar vectors. The paper states, the procedure of constructing tensors using polar vector components is ``very tedious for [tensor rank] $\ge 2$.'' And while geometrically elegant, the  use of Wigner coefficients also hinders the necessary constructions for the purpose of numerical computation. 

In our work, we find a straightforward method for making the construction iteratively using polar vectors. Our formulae involve explicit algebraic relations, and make no use of Wigner coefficients. The $\mathcal{Q}_{\ell}$ coefficients that connect the angular and radial dependencies also manifestly satisfy orthogonality relations without any kind of spatial integration. This approach makes the pseudospectral approach very straightforward. First we project out the angular polar coordinate vectors, then we perform a series of individual spin-weighted spherical harmonic transforms.  Then we project out a series of radially dependent functions in different regularity classes. Sparse differential operations can act on regularity classes, and the whole thing can go back to grid space for nonlinear multiplications. 

In \S4 we continue our development of numerical methods that work well given the geometric properties of spheres in three dimensions. For this section we use Jacobi polynomials as basis functions for $\Reg{}$. This follows almost the same pattern as V16. It happens that the operators mapping between regularity spaces have very natural action on classes of properly rescaled Jacobi polynomials. We find that it is tremendously helpful to make very specific choices regarding the Jacobi parameters. 

In \S5 we discuss how to use the results of \S2--4 to construct linear-algebraic systems for solving common PDEs. A main theme of V16 was that letting the Jacobi parameters adjust allowed for representing differential operators with maximally sparse matrices. This allows for the simple and well-conditioned construction of numerical equations from arbitrary formulations of physical equations. This not only allows for highly accurate results, but also automates the process in a way that helps reduce human error. We also show how our basis sets map to other past works. Some authors use rescaled Jacobi polynomials in one form or another, but do not allow the parameters to change as derivatives act \cite{boyd2011comparing,
livermore_jones_worland_2007,
matsushima1995spectral,
sakai2009application,
zernike}. 

In \S6 we provide conclusions for this Part-I of our two-part series. The companion paper, Part-II, picks up with a series of numerical examples showing the effectiveness of our mathematical methods. We do this with a series of numerical eigenvalue problems and full-scale nonlinear computations of numerical benchmark problems. The first set of tests shows the high accuracy in well-understood unit tests of different mathematical aspects. The second set of tests shows that everything can be combined relatively easily to make numerical solvers that might once have taken several person-years to create. We show that we can use our mathematical framework to create automated numerical solvers within the Dedalus framework.

\section{Tensor calculus on the unit 2-sphere}

In this section, we assume $r=1$. For the three-dimensional part of the paper (\S3), we assume $0 \le r \le 1$. Restricting to the unit sphere in this section eliminates pesky $1/r$ factors in the various formulae. 

\subsection{The connection}

The foremost difficulty with computations in spherical geometry is that the basis vectors change from point to point. The gradient on the 2-sphere is
\Beq
\D{} \ = \  \e{\theta}\D{\theta} + \e{\phi}\D{\phi},
\Eeq
where $\e{\theta},\e{\phi}$ represent orthonormal coordinate unit vectors on the surface of the sphere. For scalar-valued functions
\Beq
\label{grad-components}
\D{\theta} f(\theta, \phi) \ = \ \frac{\partial}{\partial \theta}f(\theta, \phi), \quad \D{\phi} f(\theta, \phi) \ = \ \frac{1}{\sin \theta}\frac{\partial}{\partial \phi}f(\theta, \phi).
\Eeq

How the covariant derivative acts on a vector sets it apart from a re-scaled set of partial derivatives. The affine connection in the coordinate basis is
\Beq
\nabla \e{\theta} \ = \ \cot \theta\, \e{\phi} \e{\phi} ,\quad \nabla \e{\phi} \ = \ - \cot \theta\, \e{\phi} \e{\theta} .
\Eeq
The $\cot \theta$ factors reflect the singular nature of the basis vectors at the poles.
We can define an alternative spinor basis that simplifies the affine connection,
\Beq
e_{\pm} = \frac{1}{\sqrt{2}}(e_{\theta} \mp i e_{\phi}).
\Eeq
We refer to Appendix A for a number of different useful properties of this basis. Most important, the $e_{\pm}$ basis diagonalises the action of the gradient such that
\Beq
\nabla e_{\pm} \ = \ \pm  i \cot \theta e_{\phi} e_{\pm}  \ = \ \pm \frac{\cot \theta}{\sqrt{2}}( e_{-} - e_{+}) e_{\pm} 
\label{spin-connection} ,
\Eeq
The gradient in the new basis becomes,
\Beq
\D{} \ \equiv \  \e{\mu}\D{\mu} \ = \  \e{+}\D{+} + \e{-}\D{-}, \label{grad-mu-def}
\Eeq
where
\Beq
\D{\pm} \ = \ \frac{1}{\sqrt{2}} \left( \D{\theta} \pm i \D{\phi} \right).
\Eeq
We use Einstein summation notation for the first part of \eq{grad-mu-def}.

For each individual component of the gradient and basis elements 
\Beq
\label{spin connection}
\D{\mu } \e{\sigma} \ = \ - \mu \sigma  \frac{\cot \theta}{\sqrt{2}} \,\e{\sigma},  
\Eeq
where $\mu,\sigma = \pm1$. So far, this analysis ignores the radial aspects, which we discuss in \S4.

We use a modified notation to denote higher-rank tensors. For a multi-index $\sigma = \{ \sigma_{1}, \ldots, \sigma_{\mathfrak{r}} \}$ with $\sigma_{i} = \pm 1$, then 
\Beq
e(\sigma) \ \equiv \ e_{\sigma_{1}}\ldots e_{\sigma_{\mathfrak{r}}}, \quad \mathrm{and} \quad \bar{\sigma} \ \equiv \ \sum_{j=1}^{\mathfrak{r}} \sigma_{j}
\label{multi-index sum}
\Eeq
See Appendix A for more details regarding this notation. Acting on an arbitrary-rank tensor component \eq{spin connection} becomes 
\Beq
\label{De}
\D{\mu } e(\sigma) \ = \ - \mu\, \overline{\sigma}\,  \frac{\cot \theta}{\sqrt{2}} e(\sigma).  
\Eeq
The presence of the spin-weight, $\bar{\sigma}$, is the only difference between the formula for single and multi-index components. 

\subsection{Spin-weighted Spherical Harmonics}

Spin-weighted spherical harmonics generalise traditional spherical harmonics. Both sets of functions result from classical Jacobi polynomials with different parameters. Written in full form, 
\Beq
Y_{\ell,m}^{s}(\theta,\phi) =  N_{\ell,m}^{s} \sin^{|m+s|}\left(\tfrac{\theta}{2}\right)  \cos^{|m-s|} \left(\tfrac{\theta}{2}\right) P_{\ell-\ell_{0}}^{(|m+s|,|m-s|)}(\cos \theta) e^{i m \phi}
\Eeq
where
\Beq
N_{\ell,m}^{s} \ = \ (-1)^{\max(m,-s)}\sqrt{
\tfrac{(2\ell+1) (\ell+ \ell_{0}  )! (\ell-  \ell_{0}  )! }{2(\ell+  \ell_{1} )!(\ell- \ell_{1} )!}},\\ \ell_{0} = \max (|m|, |s|), \quad \ell_{1} = \min (|m|, |s|).
\Eeq
The function $P_{n}^{(a,b)}(z)$ represents a Jacobi polynomial of degree $n=\ell-\ell_{0}$ and parameters $a,b = |m \pm s|$. 
We use the recursion properties of Jacobi polynomials to construct operators acting on this spherical basis.

Traditional (scalar-valued) spherical harmonics correspond to $s=0$. The covariant derivative shifts the spin-weight, $s \to s \pm 1$. Specifically, 
\Beq
\label{DY}
\D{\mu} Y^{s}_{\ell,m} -  \mu s\,\frac{\cot \theta}{\sqrt{2}} Y^{s}_{\ell,m} \ = \  k^{\mu}_{\ell,s} Y^{\mu+s}_{\ell,m}   ,
\Eeq
where
\Beq
\label{k-wavenumbers}
k^{\mu}_{\ell,s} \ \equiv \ -\mu \sqrt{\frac{(\ell- \mu s)(\ell+ \mu s+1)}{2}}.
\Eeq
The $\cot \theta $ term in \Eq{DY} implies that the covariant derivative of the spin-weighted harmonics can become singular in some cases. 
However, putting \eq{DY} together with \eq{De} gives 
\Beq
\D{\mu} \left(Y^{\bar{\sigma}}_{\ell,m}e(\sigma) \right) \ = \  k^{\mu}_{\ell,\bar{\sigma}}  Y^{\mu+\bar{\sigma}}_{\ell,m}e(\sigma),
\Eeq
where $\mu, \sigma_{i} =\pm1$. The combination of the spin basis, and spin-weighted spherical harmonics yields a very simple action for the gradient operator. For the case $\bar{\sigma}\ne0$, both $Y^{\bar{\sigma}}_{\ell,m}$ and $e(\sigma)$ each individually produce singular gradients at the poles; the product of the two remains regular everywhere.   

For an arbitrary tensor field
\Beq
\mathrm{T} = \sum_{\ell,m,\sigma} A^{\sigma}_{\ell,m} Y^{\bar{\sigma}}_{\ell,m}e(\sigma).
\Eeq
Now,
\Beq
\nabla \mathrm{T} = \sum_{\mu\sigma,\ell,m}  k^{\mu}_{\ell,\bar{\sigma}}\, A^{\sigma}_{\ell,m} Y^{\overline{\mu\sigma} }_{\ell,m }\,e(\mu \sigma),
\Eeq
where $\mu \sigma = \{\mu, \sigma_{1}, \ldots, \sigma_{\mathfrak{r}}\}$ and $\overline{\mu\sigma} = \mu + \bar{\sigma}$; see Appendix A.
The new spherical harmonic coefficients of $\nabla \mathrm{T}$ are simple multiples of the original coefficients for $\mathrm{T}$. For the spectral coefficients of a general multi-index tensor field
\Beq
\nabla_{\mu} \ : \  A^{\sigma}_{\ell,m} \ \to \ k^{\mu}_{\ell,\bar{\sigma}}\, A^{\sigma}_{\ell,m}
\Eeq
This is one of the main points of the paper:
\begin{itemize}
\item
\textsl{Using a combination of spin-weighted spherical harmonics, and spinor basis vectors, differentiation of the sphere becomes just like Fourier series; \ie purely diagonal wavenumber multiplication.}
\item
\textsl{The major difference in this case is that Fourier series uses the same basis to represent scalar functions and the components of vectors and tensors. The sphere requires a hierarchy of bases. The one-dimensional circle also uses a small hierarchy of bases, \ie $\sin$/$\cos$.}
\end{itemize}

\subsection{2-sphere Laplacian}

The intrinsic Laplacian on the sphere results from contracting two gradient operations  
\Beq
\label{rough Laplacian}
\D{}\cdot \D{}  \left(Y^{\bar{\sigma}}_{\ell,m}e(\sigma) \right)  \ = \  \left(k^{-}_{\ell,\bar{\sigma}+1} k^{+}_{\ell,\bar{\sigma}} + k^{+}_{\ell,\bar{\sigma}-1} k^{-}_{\ell,\bar{\sigma}} \right)Y^{\bar{\sigma}}_{\ell,m}e(\sigma),
\Eeq
where (for a single value of spin-weight, $s$), 
\Beq
k^{-}_{\ell,s+1} k^{+}_{\ell,s} + k^{+}_{\ell,s-1} k^{-}_{\ell,s}  \ = \ - \ell(\ell+1) + s^{2}.
\Eeq
Also note that
\Beq
k^{-}_{\ell,s+1} k^{+}_{\ell,s} - k^{+}_{\ell,s-1} k^{-}_{\ell,s}  \ = \ s.
\Eeq
Therefore the commutator of derivatives is:
\Beq
\left[\D{-}, \D{+}\right] Y^{\bar{\sigma}}_{\ell,m}e(\sigma) \ = \ \bar{\sigma}\, Y^{\bar{\sigma}}_{\ell,m}e(\sigma).
\Eeq
\Eq{rough Laplacian} gives a slightly different formula than results from taking the three-dimensional Laplacian and restricting it to the surface of the 2-sphere. In this case, additional terms result from contracting in the third dimension. \Eq{rough Laplacian} defines what is often called the \textsl{rough Laplacian}, or \textsl{connection Laplacian}. In a curved geometry, it is possible to define several linear, second-order,  elliptic differential operators with a reasonable claim to the title of Laplacian. The Weitzenb\"{o}ck identity ensures that any two such Laplacians differ by a scalar curvature term at most; \ie a term with no derivatives; see \eg \cite{Homma}.  In the case of the restricted three-dimensional Laplacian, 
\Beq
\label{Other-Laplacian}
 \D{3D}\cdot \D{3D}  \left(Y^{\bar{\sigma}}_{\ell,m}e(\sigma) \right) \big|_{S^{2}}  \ = \  -\left(\ell(\ell+1) - \bar{\sigma}^{2} + | \sigma| \right)Y^{\bar{\sigma}}_{\ell,m}e(\sigma),
\Eeq
where $| \sigma| = \mathfrak{r}$ gives the tensor rank of $e(\sigma)$. In the case of the simple vectors $e_{\pm}$, we have that $\bar{\sigma}^{2} = | \sigma | = 1$. In applications, the most appropriate Laplacian depends on the details of the underlying physics.

\subsection{Multiplication by trigonometric functions}

Numerical computations for PDEs obviously require derivative operations. Some calculations also require multiplication by non-constant coefficients. On the surface of the 2-sphere, the action of $\sin \theta$ and $\cos \theta$ allows for multiplication by arbitrary well-behaved functions. Treating $\sin\theta$ and $\cos\theta$ as operators acting between bases of spin-weighted spherical harmonics allows for non-constant coefficient terms (\eg Coriolis and shear) on the 2-sphere in time-evolution and/or eigenvalue problems.

The multiplication of non-constant coefficients act as operators in the entire basis set of spin-weighted spherical harmonics. Defining the row vector, 
\Beq
{Y}_{m}^{s} \ = \ \left[\,Y_{\ell_{0},m}^{s}, \, Y_{\ell_{0}+1,m}^{s} \, \ldots \, \right]
\Eeq
we have
\Beq
\cos (\theta)\, {Y}_{m}^{s} \ = \ {Y}_{m}^{s} \, \mathcal{C}_{s}, \quad
\sin (\theta)\, {Y}_{m}^{s} \ = \ {Y}_{m}^{s\pm1} \, \mathcal{S}^{\pm}_{s}.
\Eeq
The operators $\mathcal{C}_{s}$ and $\mathcal{S}_{s}$ have the property 
\Beq
\mathcal{C}_{s\pm1} \mathcal{S}_{s}^{\pm} \ = \ \mathcal{S}_{s}^{\pm} \mathcal{C}_{s}, \quad \mathcal{S}_{s+1}^{-}\mathcal{S}_{s}^{+} = \mathcal{S}_{s-1}^{+}\mathcal{S}_{s}^{-},\quad \mathcal{C}_{s} \mathcal{C}_{s} + \mathcal{S}_{s-1}^{+}\mathcal{S}_{s}^{-} = \mathcal{I}_{s}.
\Eeq
where $\mathcal{I}_{s}$ represents the identity acting on the spin-weighted bases. 
\Beq
(\mathcal{S}_{s}^{+})^{\top} (\mathcal{S}_{s+1}^{-})^{\top} = (\mathcal{S}_{s}^{-})^{\top} (\mathcal{S}_{s-1}^{+})^{\top}
\Eeq
The sine multiplication operators are related via the transpose, while the cosine operators are symmetric. That is, 
\Beq
\mathcal{S}_{s}^{-} \ = \  (\mathcal{S}_{s-1}^{+})^{\top} \quad \mathrm{and} \quad \mathcal{C}_{s} \ = \ \mathcal{C}_{s}^{\top}
\Eeq
The $\mathcal{C}_{s}$ is the Jacobi matrix needed for the construction of the spin-weighted spherical harmonics themselves. 

Explicitly, in terms of each basis element, 
\Beq
\label{cos times}
\cos \theta\, Y^{s}_{\ell,m} & = & \beta_{\ell+1} \, Y^{s}_{\ell+1,m} + \alpha_{\ell}\, Y^{s}_{\ell,m} + \beta_{\ell}\,  Y^{s}_{\ell-1,m}\\
\sin \theta\, Y^{s}_{\ell,m} &=& \kappa_{\ell+1,-s-1} \, Y^{s+1}_{\ell+1,m} + \rho_{\ell,s} \, Y^{s+1}_{\ell,m}  - \kappa_{\ell,s}\,  Y^{s+1}_{\ell-1,m} \\ 
 \sin \theta\, Y^{s}_{\ell,m} &=& -\kappa_{\ell+1,s-1} \, Y^{s-1}_{\ell+1,m} + \rho_{\ell,s-1} \, Y^{s-1}_{\ell,m}  + \kappa_{\ell,-s}\,  Y^{s-1}_{\ell-1,m}
\label{sin times}
\Eeq
where
\Beq
\alpha_{\ell} \ = \ -\frac{m s}{\ell(\ell+1)}, & \beta_{\ell} \ = \ \frac{1}{\ell }\sqrt{\tfrac{(\ell -m)
   (\ell +m) (\ell -s) (\ell+s
   )}{(2 \ell -1) (2 \ell
   +1)}} \label{Cos Coefficients} \\
 \rho_{\ell,s} \ = \ m\tfrac{ \sqrt{(\ell -s) (\ell + s
   +1)}}{\ell  (\ell +1)},  &  \kappa_{\ell,s} \ = \ \frac{1}{\ell}\sqrt{\tfrac{(\ell -m)
   (\ell +m) (\ell-s -1) (\ell
   -s)}{(2 \ell -1) (2 \ell
   +1)}} \label{Sin Coefficients}
\Eeq
The special case $\ell=m=0$ in \eqs{Cos Coefficients}{Sin Coefficients} is handled by setting $m=0$ and cancelling terms before taking $\ell=0$. Cosine multiplication preserves the spin-weight, while sine multiplication can either raise or lower the spin-weight depending on the situation. We include an example using the cosine operators in Paper-II.

\section{Tensor calculus of the 3-ball}

The main purpose of this paper is to use, decompose and organise all scalar, vector, and tensor computations in the three-dimensional ball into a mathematical structure that is efficient for numerical computations. The previous section introduced how to do this for the angular components. This section lays out a similar framework for the radial direction.

\subsection{Motivation: Scalar functions}

For a scalar field, we can decompose the radial and angular coordinates into a function of the form
\Beq
\label{f-example}
f(r,\theta,\phi) \ = \ F (r) Y^{0}_{\ell,m}(\theta,\phi).
\Eeq
The most general situation would involve a sum over $m$, and $\ell$, but \eq{f-example} is enough to illustrate the main points. The Laplacian acts in the radial direction such that 
\Beq
\nabla^{2} f(r,\theta,\phi) \ = \ \left[ \frac{1}{r^{2}}\frac{d}{dr}\left(r^{2} \frac{d}{dr}F(r)\right)  - \frac{\ell(\ell+1)F(r)}{r^{2}} \right] Y^{0}_{\ell,m}
\Eeq
The second-order radial operator decomposes into the product of two $\ell$-dependent first-order operators, $D_{\ell}^{\pm}$, such that
\Beq
\label{Laplacian=DD}
\left[ \frac{1}{r^{2}}\frac{d}{dr}\left(r^{2} \frac{d}{dr}F(r)\right)  - \frac{\ell(\ell+1)F(r)}{r^{2}} \right]  \ = \ D_{\ell+1}^{-}D_{\ell}^{+}F(r) \ = \ D_{\ell-1}^{+}D_{\ell}^{-}F(r)
\Eeq
where
\Beq
\label{D+,D-}
D_{\ell}^{+} \ \equiv \ \frac{d}{dr} - \frac{\ell}{r}, \quad D_{\ell}^{-} \ \equiv \ \frac{d}{dr} + \frac{\ell+1}{r}.
\Eeq
\Eqs{Laplacian=DD}{D+,D-} are well known in many contexts.  For example, the spherical Bessel functions $j_{\ell}(r)$ have
\Beq
j_{\ell}(r) \in \Reg{}\quad \mathrm{and} \quad D_{\ell}^{\pm} j_{\ell}(r) \ = \ \mp j_{\ell\pm1}(r). \label{Bessel-pm}
\Eeq
These operators are the spherical versions of the $m$-dependent operators we discussed for the 2D disk. In the disk, we were able to use the two factors of the Laplacian to create components of the gradient operator. But the sphere presents us with a conundrum: in the sphere there are two factors of the Laplacian, but three components of the gradient. Furthermore, the angular gradient on spherical harmonic functions does not produce the effect of multiplying by $\ell$ or $\ell+1$. This is in contrast to the disk where the angular derivative on exponential functions does result in the multiplication by $m$. How can we use the operators in \eq{D+,D-} in a general tensor framework? 

The main reason that these operators in (\ref{D+,D-}) are useful is the fact that
\Beq
D_{\ell}^{\pm} \ : \ \Reg{} \ \to \ \mathrm{Reg}(\ell\pm1)
\Eeq
We will see this more explicitly in \S4 when considering computational strategies.  The $D_{\ell}^{\pm}$ operations appear individually for vector and tensors computations, but getting to the point where these operators appear takes some effort compared to scalars.  

\subsection{Scalars and Vectors}

\subsubsection{The connection}

In three dimensions we have an additional basis vector $e_{0} \equiv e_{r}$. The gradient formulae for the spin basis now acquires additional terms arising from the radial element,
\Beq
\label{3D-connection}
\D{\mu}\e{\sigma} \ = \ - \mu \sigma \frac{\cot\theta}{\sqrt{2} r}\e{\sigma} - \frac{\delta_{\mu,\sigma}}{r}\e{0} + \frac{\delta_{\sigma,0}}{r}\e{-\mu}.
\Eeq
\Eq{3D-connection} reduces to \eq{spin connection} when restricted to the surface of the 2-sphere (\eg $r=1$). 

A difficulty exists with three dimensions that is not present with the  2-sphere: there is no single rotation (complex-valued or otherwise) that diagonalises \eq{3D-connection} in a manner similar to \eq{spin connection}. In the 2-sphere, the $e_{\pm}$ basis diagonalises the connection coefficients. This is seen in \eq{3D-connection} where the first term on the right-hand side is proportional to $e_{\sigma}$, which is the input on the left-hand side. To make sense of the other terms, it is helpful to consider the action of $\D{\mu}$ on all three basis vectors simultaneously. Therefore, define
\Beq
\mathcal{E} \ \equiv \ \left[
\begin{array}{c}
 e_{-} \\ e_{0} \\  e_{+} 
\end{array}
\right].
\Eeq
The first thing we can say about \eq{3D-connection} is that
\Beq
\D{0}\mathcal{E} \ = \ 0. \label{grad 0 basis}
\Eeq
The other components are equivalent to 
\Beq
\D{\pm} \mathcal{E} \ = \ \mp    \frac{1}{\sqrt{2} r} \left(\cot \theta\, \mathcal{J}_{0} + \mathcal{J}_{\pm} \right) \mathcal{E}  , \label{grad +- basis}
\Eeq
where
\Beq
\mathcal{J}_{0} \ \equiv \ \left[
\begin{array}{ccc}
 -1 & 0 & 0 \\
 0 & 0 & 0 \\
 0 & 0 & 1 \\
\end{array}
\right]
\Eeq
and
\Beq
\label{sp(2,r)}
 \mathcal{J}_{-} \ \equiv \ \left[
\begin{array}{ccc}
 0 & -\sqrt{2} & 0 \\
 0 & 0 & \sqrt{2} \\
 0 & 0 & 0 \\
\end{array}
\right], \quad  \mathcal{J}_{+} \ \equiv \ \left[
\begin{array}{ccc}
 0 & 0 & 0 \\
 -\sqrt{2} & 0 & 0 \\
 0 & \sqrt{2} & 0 \\
\end{array}
\right].
\Eeq
We see that $\mathcal{J}_{0}$ is diagonal.  This is what we mean when we say that the $e_{\pm}$ diagonalises the connection in the angular directions. Also, clearly $ \mathcal{J}_{\pm}^{\dag} =  \mathcal{J}_{\mp}$. The question is, how do we know there is not a transformation that diagonalises all three matrices? This would simplify tensor calculus in three dimensions considerably. We can simultaneously diagonalises two matrices if and only if they commute. Checking the commutation relations, 
\Beq
\left[ \mathcal{J}_{0}, \mathcal{J}_{\pm} \right] \ = \ \pm \mathcal{J}_{\pm}, \quad \left[ \mathcal{J}_{+}, \mathcal{J}_{-} \right] \ = \ 2 \mathcal{J}_{0}.
\Eeq
No only do the matrices not commute, but they form the Lie algebra $\mathfrak{sl}(2) \simeq \mathfrak{sp}(2)   $.
Transformations do exist that simplify computations, but they require additional mathematical ingredients. We need to work in spectral space to find what we are looking for.

\subsubsection{Motivation: Vector Laplacian.}

We can represent 3D vector fields in terms of spin-weighted spherical harmonics in a method similar to 2D vector fields and/or scalar fields. For example,  
\Beq
\label{v-example}
v(r,\theta,\phi) = V^{-}(r) Y^{-}_{\ell,m}(\theta,\phi) e_{-}+ V^{0}(r) Y^{0}_{\ell,m}(\theta,\phi) e_{0} + V^{+}(r) Y^{+}_{\ell,m}(\theta,\phi) e_{+}.
\Eeq
Like \eq{f-example} the most general situation could involve a sum over $m$, and $\ell$, but this is enough to illustrate the main points.
The Laplacian now acts in a more complicated way than for scalars,
\Beq
& e_{-}^{*} \cdot \nabla^{2} v  \ = \ \left[\frac{1}{r^{2}}\frac{d}{dr}(r^{2} \frac{d}{dr}V^{-})  - \frac{\ell(\ell+1) V^{-}}{r^{2}} + \frac{\sqrt{2\ell(\ell+1)}V^{0}}{r^{2}}\right] Y_{\ell,m}^{-},  \\
&  e_{0}^{*} \cdot \nabla^{2} v \ = \  \left[\frac{1}{r^{2}}\frac{d}{dr}(r^{2} \frac{d}{dr}V^{0}) - \frac{(\ell(\ell+1)+2) V^{0}}{r^{2}} + \frac{\sqrt{2 \ell (\ell+1)}( V^{-} - V^{+} )}{r^{2}}\right]  Y_{\ell,m}^{0}, \\
& e_{+}^{*} \cdot \nabla^{2} v \ = \   \left[\frac{1}{r^{2}}\frac{d}{dr}(r^{2} \frac{d}{dr}V^{+}) - \frac{\sqrt{2\ell(\ell+1)}V^{0}}{r^{2}} - \frac{\ell(\ell+1) V^{+}}{r^{2}}   \right] Y_{\ell,m}^{+}.
\Eeq
Restricting this to the surface of the 2-sphere coincides with \eq{Other-Laplacian}. 
The vector Laplacian couples the different spin-weight components such that
\Beq
\nabla^{2} v \ = \ \left[\frac{1}{r^{2}}\frac{d}{dr}(r^{2} \frac{d}{dr}V^{\mu}) - \frac{1}{r^2}\Lambda_{\nu}^{\mu} V^{\nu}\right] Y^{\bar{\mu}}_{\ell,m}e_{\mu},
\Eeq
with the $\ell$-dependent symmetric matrix 
\Beq
\Lambda \ \equiv \ \left[
\begin{array}{ccc}
   \ell
   (\ell+1) &-\sqrt{2\ell (\ell+1)} & 0 \\
     -\sqrt{2\ell
   (\ell+1)} & \ell(\ell+1)+2 &  \sqrt{2\ell
   (\ell+1)} \\
 0 &  \sqrt{2\ell (\ell+1)} &
   \ell (\ell+1) \\
\end{array}
\right].
\Eeq
The eigenvalues of $\Lambda$ are
\Beq
\lambda_{\ell-1}, \quad \lambda_{\ell}, \quad \lambda_{\ell+1}\, \quad \mathrm{where} \quad \lambda_{\ell} \equiv \ell(\ell+1)
\Eeq
We see clearly that the Laplacian acts on vectors as if they were composed of three scalars each with $\ell-1$, $\ell$, $\ell+1$. This is what we mean when we say that the components of vectors are in a direct sum of regularity classes. The eigenvectors of $\Lambda$ are orthogonal. The following matrix diagonalises the Laplacian in $\ell$, 
\Beq
\label{Q-matrix}
\mathcal{Q}_{\ell}  \ \equiv \ \frac{1}{\sqrt{2}\sqrt{2 \ell+1}}\left[
\begin{array}{ccc}
 \sqrt{\ell+1} & \sqrt{2 \ell+1} & -\sqrt{\ell} \\
  \sqrt{2\ell} & 0 &  \sqrt{2(\ell+1)}
   \\
 -\sqrt{\ell+1} & \sqrt{2 \ell+1} & \sqrt{\ell} \\
\end{array}
\right].
\Eeq
That is,   
\Beq
\mathcal{Q}_{\ell} ^{\top} .\, \Lambda  .\, \mathcal{Q}_{\ell}  \ = \ \left[
\begin{array}{ccc}
 \lambda_{\ell-1} & 0 & 0 \\
 0 & \lambda_{\ell} & 0 \\
 0 & 0 & \lambda_{\ell+1} \\
\end{array}
\right].
\Eeq
We use the notation $\mathfrak{B}(s)$ of Eastwood \& Tod (1982) to denote spaces of coherent $\textit{spin-weight} = s$. That is, $Y_{\ell,m}^{s}$ forms an orthonormal basis for $\mathfrak{B}(s)$. The matrix has the property that,
\Beq
\mathcal{Q}_{\ell} \quad : \quad \Reg{-1}\oplus \Reg{}\oplus \Reg{+1} \quad \to \quad \mathfrak{B}(-) \oplus \mathfrak{B}(0) \oplus \mathfrak{B}(+), 
\Eeq
and
\Beq
\mathcal{Q}^{\top}_{\ell} \quad : \quad \mathfrak{B}(-) \oplus \mathfrak{B}(0) \oplus \mathfrak{B}(+) \quad \to \quad \Reg{-1}\oplus \Reg{}\oplus \Reg{+1},
\Eeq 
where the notation $\oplus$ in this case represents the \textit{direct sum} of vector spaces\footnote{It is standard to use the $\oplus$ symbol \textit{both} for the direct sum, and for the Kronecker sum. We use Kronecker sum in \S3.3.2.}. 
We therefore arrive at another main point of this paper. 
\begin{itemize}
\item
\textsl{Finding simple and sparse operators in the radial direction requires an orthogonal rotation on the vector components. This rotation is spherical-harmonic-degree ($\ell$) dependent.}
\end{itemize}
Because $\mathcal{Q}_{\ell}$ depends on $\ell$ its action must occur after a spherical-harmonics transformation. 

In defining the matrix $\mathcal{Q}_{\ell}$, it is useful to define the following notation 
\Beq
\label{grad-unit}
\xi_{\ell}^{-} \ \equiv  \ \sqrt{ \frac{\ell}{2\ell+1}},\quad  \xi_{\ell}^{+} \ \equiv  \ \sqrt{ \frac{\ell+1}{2\ell+1}},\quad \mathrm{and} \quad \xi_{\ell}^{0} \ \equiv  \ 0.
\Eeq
These factors appear in many of the following formulae. They have the useful property that
\Beq
\left( \xi_{\ell}^{-} \right)^2 + \left( \xi_{\ell}^{+} \right)^2 \ = \ 1
\Eeq
With this notation, \eq{Q-matrix} becomes
\Beq
\label{Q-xi-matrix}
\mathcal{Q}_{\ell}  \ = \ \left[
\begin{array}{ccc}
 \tfrac{1}{\sqrt{2}}\xi_{\ell}^{+} & \tfrac{1}{\sqrt{2}} & -\tfrac{1}{\sqrt{2}}\xi_{\ell}^{-} \\
  \xi_{\ell}^{-} & 0 &  \xi_{\ell}^{+}
   \\
 -\tfrac{1}{\sqrt{2}}\xi_{\ell}^{+} & \tfrac{1}{\sqrt{2}} & \tfrac{1}{\sqrt{2}}\xi_{\ell}^{-}  \\
\end{array}
\right].
\Eeq

\textsl{Notation} --- In working with spin-weight degrees of freedom, and regularity degrees of freedom, we use particular notation to denote each type of index. We use Greek letters to denote spin indices, \eg $\mu$, $\nu$, $\sigma$. We use non-italics Roman letters toward the beginning of the alphabet to denote regularity class relative to $\ell$, \eg $\mathrm{a}$,  $\mathrm{b}$, $\mathrm{c}$. Also, we use a function notation to represent the indices of the orthogonal matrices, \eg $\mathcal{Q}_{\ell}(\sigma, \mathrm{a})$ represents the $\sigma$-th, $\mathrm{a}$-th element of $\mathcal{Q}_{\ell}$. We do this to avoid too many subscripts and to generalise the notation more easily to multi-index notation in higher-ranks tensors.  When we use multi-index notation, an over-bar, \eg $\bar{\sigma}$, represents a sum over multi-indices, as in equation~(\ref{multi-index sum}) and Appendix A.

\subsubsection{Vector calculus.}

Now define a scalar field $f$ and two seperate vector fields $v$ and $u$:
\Beq
f &=& F(r) Y^{0}_{\ell,m}(\theta,\phi), \\
 v &=& \sum_{\sigma,\mathrm{a}}\mathcal{Q}_{\ell}(\sigma,\mathrm{a})V^{\mathrm{a}}(r)  Y^{\sigma}_{\ell,m}(\theta,\phi) e_{\sigma}, \\  
 u & = & \sum_{\sigma,\mathrm{a}}\mathcal{Q}_{\ell}(\sigma,\mathrm{a})U^{\mathrm{a}}(r)  Y^{\sigma}_{\ell,m}(\theta,\phi) e_{\sigma},
\Eeq
where we again suppress more-general summations over $\ell$ and $m$.  Our primary vector calculus operators are then:

\textsl{Grad} --- If $v = \nabla f $, then
\Beq
\label{grad -}
&& V^{-} \ = \    \xi_{\ell}^{-} D_{\ell}^{-} F, \\ 
\label{grad 0}
&& V^{0} \ = \    0, \\
\label{grad +}
&& V^{+} \ = \    \xi_{\ell}^{+}  D_{\ell}^{+}F,
\Eeq

\textsl{Div} --- If $f = \nabla \cdot v$ then
\Beq
\label{div}
F = \xi_{\ell}^{-} D^{+}_{\ell-1}V^{-} + \xi_{\ell}^{+} D^{-}_{\ell+1}V^{+}.
\Eeq

\textsl{Curl} --- If $u = \nabla \times v$, then
\Beq
\label{curl -}
&& U^{-} \ = \ -i \xi_{\ell}^{+} D^{-}_{\ell} V^{0}
\\
\label{curl 0}
&& U^{0} \ =\  i \xi_{\ell}^{-} D^{-}_{\ell+1}V^{+} -i \xi_{\ell}^{+} D^{+}_{\ell-1}V^{-}
\\
\label{curl +}
&&U^{+} \ = \  i \xi_{\ell}^{-} D^{+}_{\ell}V^{0}
\Eeq

\textsl{Laplacian} --- If $u = \nabla^{2} v$, then different iterated combinations of \eqss{grad -}{curl +} imply that $u =  \nabla (\nabla \cdot v ) - \nabla \times (\nabla \times v)$. And also,
\Beq
U^{\mathrm{a}} \ = \ D^{-}_{\ell+\mathrm{a}+1}  D^{+}_{\ell+\mathrm{a}}V^{\mathrm{a}} \ = \ D^{+}_{\ell+\mathrm{a}-1}  D^{-}_{\ell+\mathrm{a}}V^{\mathrm{a}}
\Eeq
which is exactly the same form as the diagonalised vector Laplacian from \S{3.2.2}.

The divergence of the curl, and the curl of the gradient vanish. 

\subsection{General Tensor Calculus}

We want to find a way to operate on general tensors in a similar way we can for vectors and scalars. This is important in many physical applications. Stress tensors appear naturally both in solid and fluid mechanics, and magnetohydrodynamics. General relativity requires 4th-order tensor fields. The polarisation field arises naturally as a  2nd-rank tensor in electro-magnetics. 

We first try to derive a simple set of transformations from the Laplacian in the same way we did for vector fields. This turns out to be informative, but does not completely determine everything we need. A systematic approach requires iteratively apply the gradient operator to a tensor field and cleaning up the result. 

\subsubsection{Motivation: Tensor Laplacian}

For a single $\ell,m$ component of a 2-tensor
\Beq
\label{T-example}
T(r,\theta,\phi) &= & \sum_{\sigma=-1}^{+1} \sum_{\tau=-1}^{+1} T^{\sigma\tau}(r) Y^{\sigma+\tau}_{\ell,m}(\theta,\phi) e_{\sigma}e_{\tau}
\Eeq
After separating out the angular components, the 2-tensor Laplacian becomes
\Beq
\nabla \cdot \nabla\, T \ \to \ \frac{1}{r^{2}}\frac{d}{dr}( r^{2} \frac{d}{dr } T^{\sigma\tau}) - \frac{1}{r^{2}}\Lambda^{\sigma\tau}_{\mu\nu}\, T^{\mu\nu}
\Eeq
In this case, $\Lambda$ is a 4th-rank operator (indexed with the multi-index $\sigma\tau,\mu \nu$) that acts like a matrix on the 2nd-rank tensors. The operator is symmetric in the sense that $\Lambda^{\sigma\tau}_{\mu\nu} = \Lambda^{\mu\nu}_{\sigma\tau}$. Similar to matrices,  $\Lambda$  has eigenvalues in the sense that 
\Beq
\label{Tensor-eigen}
\Lambda^{\sigma\tau}_{\mu\nu}\, T^{\mu\nu} \ = \ \lambda_{\ell+\mathrm{a}}\, T^{\sigma\tau},
\Eeq 
where $ -2 \le \mathrm{a} \le 2$.
\Eq{Tensor-eigen} represents $3^2 = 9$ seperate linear equations. The eigenvalues fall into degenerate classes such that
\Beq
\lambda_{\ell+\mathrm{a}} \ \in \ \left\{ \lambda_{\ell-2},\ \lambda_{\ell-1},\ \lambda_{\ell-1},\ \lambda_{\ell},\ \lambda_{\ell},\ \lambda_{\ell},\  \lambda_{\ell+1},\ \lambda_{\ell+1}, \ \lambda_{\ell+2} \right\}.
\Eeq
This implies that for the spherical-harmonic coefficients of rank-2 tensors, 
\Beq
T_{\ell}(r) \ \in \ \Reg{-2}\oplus \Reg{-1}^{2} \oplus \Reg{}^{3}\oplus\Reg{+1}^{2} \oplus \Reg{+2}
\Eeq
For rank-$\mathfrak{r}$ tensors, the number of eigenvalue degeneracies for each spin-weight, $-\mathfrak{r} \le \mathrm{a} \le \mathfrak{r}$,  follows a generalised Pascal's triangle,
\Beq
d(\mathrm{a},\mathfrak{r}+1) =d(\mathrm{a}-1,\mathfrak{r}) + d(\mathrm{a},\mathfrak{r}) + d(\mathrm{a}+1,\mathfrak{r}).
\Eeq
Unlike the traditional Pascal triangle, there are no simple closed-form expressions for this recursion relation. The coefficients do satisfy a trinomial expansion
\Beq
(1+x+x^2)^{\mathfrak{r}} \ = \ \sum_{q=0}^{2\mathfrak{r}} d(q-\mathfrak{r},\mathfrak{r})\, x^{q}.
\Eeq
This is the generating function for the number of ways to select $q=\mathrm{a}+\mathfrak{r}$ objects from $\mathfrak{r}$ objects where we can select any individual object up to two times. 

In the case of vectors, the Laplacian acts on spaces of non-degenerate eigenvectors; \ie $\ell,\ell\pm 1$.  For tensors, the picture is not as clear. We want to build up higher-order transformations that serve the same purpose as the $\mathcal{Q}_{\ell}$ matrices in \eq{Q-matrix}. For this we need to define an efficient way to organise higher-order tensor indices. 

%\subsubsection{Second motivation: Vector gradient}

%
%%\Beq
%%u & = & \sum_{\tau, \mathrm{b} = -1}^{+1}e_{\tau} Y^{\tau}_{\ell,m}(\theta,\phi)  \mathcal{Q}_{\ell}(\tau,\mathrm{b})U^{\mathrm{b}}(r)  
%%\Eeq
%%
%%\Beq
%%\nabla u & = & \sum_{\tau, \mathrm{b} = -1}^{+1}   \left[  \nabla \left(e_{\tau} Y^{\tau}_{\ell,m}   \right) \mathcal{Q}_{\ell}(\tau,\mathrm{b}) U^{\mathrm{b}}(r)  + e_{0}  e_{\tau} Y^{\tau}_{\ell,m}\mathcal{Q}_{\ell}(\tau,\mathrm{b})  \frac{d U^{\mathrm{b}}(r)}{dr} \right]  
%%\Eeq
%%
%%\Beq
%% \nabla \left(e_{\tau} Y^{\tau}_{\ell,m}   \right) \ = \ e_{\pm} \nabla_{\pm} \left(e_{\tau} Y^{\tau}_{\ell,m}   \right)
%%\Eeq
%%\Beq
%%e_{\sigma} \nabla_{\sigma} \left(e_{\tau} Y^{\tau}_{\ell,m}   \right) \ = \ e_{\sigma} e_{\tau} k_{\ell,\tau}^{\sigma} Y^{\sigma+\tau}_{\ell,m} -  \sigma \frac{1}{\sqrt{2} r} \mathcal{J}_{\pm}(\sigma,\sigma')\,e_{\sigma'}
%%\Eeq

\subsubsection{Orthogonal rotations of general tensor components}

For the general case, we use a multi-index notation similar to that for the basis elements; see Appendix A.

Scalars have no associated basis vectors or indices; \ie $\sigma = \{\, \}$, and $e(\sigma) = e(\{\, \}) = 1$. Therefore define
\Beq
\mathcal{Q}_{\ell}( \{\, \}, \{\, \}) \ \equiv \ 1
\Eeq
For vectors, the multi-index consists of a single index, \ie $\sigma = \{\, \sigma_{1}\}$, with $\sigma_{1} = -1,0,+1$. A basis element is therefore $e(\sigma) = e_{\sigma_{1}}$. 

We define the higher-rank cases recursively. Assuming $\mathcal{Q}_{\ell}(\tau,\mathrm{b})$ is the orthogonal transformation for a rank-$\mathfrak{r}$ tensor, then $\mathcal{Q}_{\ell}(\sigma\tau,\mathrm{ a b})$ is the transformation for a rank-$(\mathfrak{r}+1)$ tensor. For each $\ell$, and for $\sigma,\mathrm{a} = -1,0,+1$, we define the following tensor-index recursion
\Beq
\label{Q-recursion}
\mathcal{Q}_{\ell}(\sigma\tau,\mathrm{ a b})  \ = \ 
\begin{cases} \ 
  \frac{  (\ell+\bar{\mathrm{b}}) \delta_{\sigma,0} \mathcal{Q}_{\ell}(\tau,\mathrm{b}) -  \mathcal{R}_{\ell}({\sigma\tau},\mathrm{b})}{\sqrt{(\ell+\bar{\mathrm{b}})(2(\ell+\bar{\mathrm{b}})+1)}}   & \mathrm{if}  \quad \mathrm{a} = -1 \\ \\ 
 \  \frac{\sigma  \mathcal{R}_{\ell}(\sigma\tau,\mathrm{b}) }{ \sqrt{(\ell+\bar{\mathrm{b}})(\ell+\bar{\mathrm{b}}+1)} }   & \mathrm{if} \quad  \mathrm{a} = 0 \\  \\
 \ \frac{  (\ell+\bar{\mathrm{b}}+1) \delta_{\sigma,0} \mathcal{Q}_{\ell}({\tau},{\mathrm{b}})  +  \mathcal{R}_{\ell}(\sigma\tau,\mathrm{b}) } { \sqrt{(\ell+\bar{\mathrm{b}}+1)(2(\ell+\bar{\mathrm{b}})+1)} }  & \mathrm{if}\quad \mathrm{a} = +1
\end{cases}.
\Eeq
We define the coefficients, 
\Beq
 \mathcal{R}_{\ell}(\sigma\tau,\mathrm{b})  \ \equiv \ -k_{\ell,\bar{\tau}}^{\sigma} \mathcal{Q}_{\ell}(\tau,\mathrm{b}) + \frac{\sigma}{\sqrt{2}} \sum_{\tau'}  \mathcal{J}^{\sigma}_{|\tau|}(\tau,\tau') \mathcal{Q}_{\ell}(\tau',\mathrm{b}),
\Eeq
where recursively,
\Beq
\label{tensor-J}
\mathcal{J}^{\sigma}_{\mathfrak{r}} \ = \ \mathcal{J}^{\sigma}_{\mathfrak{r}-1} \boldsymbol{\oplus}  \mathcal{J}^{\dag}_{\sigma}  \ \equiv \ \mathcal{J}^{\sigma}_{\mathfrak{r}-1} \boldsymbol{\otimes} \, \mathcal{I}_{1} \ + \  \mathcal{I}_{\mathfrak{r}-1} \boldsymbol{\otimes}\, \mathcal{J}^{\dag}_{\sigma}.
\Eeq
Recall that $\tau = \{\tau_{1}, \ldots, \tau_{\mathfrak{r}}\}$ and $\mathrm{b} = \{\mathrm{b}_{1}, \ldots, \mathrm{b}_{\mathfrak{r}}\}$ each represent multi-indices and therefore the updated multi-indices are $\sigma\tau = \{\sigma,\tau_{1}, \ldots, \tau_{\mathfrak{r}}\}$, $\mathrm{a}\mathrm{b} = \{\mathrm{a},\mathrm{b}_{1}, \ldots, \mathrm{b}_{\mathfrak{r}}\}$.  Also recall that $\bar{\mathrm{b}}$ is a sum over multi-indices, as in equation~(\ref{multi-index sum}) and Appendix~A.

In \eq{tensor-J}, 
\Beq
\mathcal{J}^{\sigma}_{1} \ \equiv \ \mathcal{J}^{\dag}_{\sigma} \ =  \  \mathcal{J}_{-\sigma}
\Eeq
is one of the 3-dimensional representation matrices of $\mathfrak{sl}(2)$ defined in \eq{sp(2,r)}. By definition $\mathcal{J}^{\sigma}_{0} = 1$. The matrices $\mathcal{I}_{\mathfrak{r}}$ represent identity matrices for $\mathrm{rank} = \mathfrak{r}$; \ie a $3^{\mathfrak{r}} \times 3^{\mathfrak{r}}$ identity matrix. The symbols $\boldsymbol{\oplus}$ and $\boldsymbol{\otimes}$ represent the Kronecker sum and product respectively\footnote{This is not to be confused with the direct sum, which we use when composing regularity spaces. }. Many modern numerical libraries enable the efficient computation of Kronecker products \cite{Kronecker}.
The derivative parameters are defined the same as for the 2-sphere, \ie \eq{k-wavenumbers}.
The matrices in \eq{tensor-J} also satisfy the $\mathfrak{sl}(2)$ Lie algebra, which is useful for proving many relations involving the $\mathcal{Q}_{\ell}$ transformations. 

\Eq{Q-recursion} defines a sequence of orthonormal (arbitrary-rank) tensors such that for each $\ell$, 
\Beq
\sum_{\mathrm{a}} \mathcal{Q}_{\ell}({\sigma},{\mathrm{a}})\mathcal{Q}_{\ell}({\sigma'},{\mathrm{a}})= \delta({\sigma,\sigma'}), \quad \sum_{\sigma} \mathcal{Q}_{\ell}({\sigma},{\mathrm{a}})\mathcal{Q}_{\ell}({\sigma},{\mathrm{a}'}) = \delta(\mathrm{a},\mathrm{a}'),
\Eeq
where $\delta({\sigma,\sigma'})$ represents a multi-index version of the Kronecker-$\delta$.

For rank-1, the individual matrix entries associated with \eq{Q-xi-matrix} are
\Beq
\label{Q-rank-1}
\mathcal{Q}_{\ell}(\sigma,\mathrm{a})  \ = \ 
\begin{cases} \ 
 \ \xi_{\ell}^{-} \delta_{\sigma,0}    -\frac{\sigma    }{\sqrt{2}} \xi_{\ell}^{+}  & \mathrm{if}  \quad \mathrm{a} = -1 \\ \\ 
 \ \frac{1 -\delta_{\sigma , 0}      }{ \sqrt{2} }   & \mathrm{if} \quad  \mathrm{a} = 0 \\  \\
 \ \xi_{\ell}^{+} \delta_{\sigma,0}    + \frac{\sigma    }{\sqrt{2}} \xi_{\ell}^{-}  & \mathrm{if}\quad \mathrm{a} = +1
\end{cases}
\Eeq
When convenient, we identify a multi-index with a single element with that element, \ie $\sigma \leftrightarrow \{\sigma\}$. \Eq{Q-rank-1} shows that \eq{Q-recursion} reduces to \eq{Q-xi-matrix} for rank-1.

The transformations in \eq{Q-recursion} do not exist if:
\Beq
\ell< \max( \sigma \bar{\tau}, - \sigma \bar{\tau} -1)   \quad \mathrm{or} \quad ( \ell= 0 \quad \mathrm{and} \quad \mathrm{a} \ < \ 1).
\Eeq
This means that some cases need exceptions for low $\ell$ values depending on the tensor rank. For example, there is only one linearly independent $\ell=0$  vector, rather than the expected 3 components, \ie $v(r) = \nabla \Phi(r)$. This means that for rank-1 and $\ell=0$, $\mathcal{Q}_{0} \, : \, \mathrm{Reg}(1) \to \mathfrak{B}(0)$ only. None of the other transformations/spaces exist. These special low-$\ell$ cases are straightforward to catch and avoid when creating the transformations. 

\Eq{Q-recursion} has the advantage of not using Wigner 3-$j$ coefficients. The formulae benefit from using explicit and elementary expressions at each stage. The formula is simple to program and stable for all values of $\ell$ and all tensor ranks, $\mathfrak{r}$. It runs efficiently within practical limits; \eg even rank-4 tensors (comprising 81 individual fields in 3D) are sufficient for nearly any realistic physical application.

\subsubsection{General Tensor Calculus}

We can define an arbitrary-rank tensor field such that
\Beq
\label{general-tensor}
\mathrm{T}(r,\theta,\phi) \ = \ \sum_{\mathrm{a},\sigma} \sum_{\ell,m}  \mathcal{Q}_{\ell}(\sigma,\mathrm{a})  A_{\ell,m}^{\mathrm{a}}(r)Y^{\bar{\sigma}}_{\ell,m}(\theta,\phi)e(\sigma)
\Eeq

\textsl{Projection} --- We need to extract the radial coefficients in order to preform any derivative operations.
Therefore,
\Beq
e( \sigma^{*} ) \cdot \mathrm{T} = \sum_{\mathrm{a}}\sum_{\ell,m}\mathcal{Q}_{\ell}(\sigma,\mathrm{a})  A_{\ell,m}^{\mathrm{a}}(r)Y^{\bar{\sigma}}_{\ell,m}(\theta,\phi).
\Eeq
See Appendix A, equation~(\ref{sigma-star}), for details on the meaning of $\sigma^{*}$.
Furthermore, 
\Beq
\int_{S^{2}}  {Y}^{\,\bar{\sigma}}_{\ell,m}(\theta,-\phi) \, e( \sigma^{*} ) \cdot \mathrm{T}(r,\theta,\phi) \,\mathrm{d} \hat{n} \ = \ \sum_{\mathrm{a}}\mathcal{Q}_{\ell}(\sigma,\mathrm{a})  A_{\ell,m}^{\mathrm{a}}(r)
\Eeq
where $\mathrm{d} \hat{n} \equiv \sin \theta\, \mathrm{d} \theta\, \mathrm{d} \phi$, and ${Y}^{\,\bar{\sigma}}_{\ell,m}(\theta,-\phi)$ implies complex conjugation. 
Lastly,
\Beq
\sum_{\sigma} \mathcal{Q}_{\ell}(\sigma,\mathrm{a})  \int_{S^{2}}  {Y}^{\,\bar{\sigma}}_{\ell,m}(\theta,-\phi) \, e( \sigma^{*} ) \cdot \mathrm{T}(r,\theta,\phi) \,\mathrm{d} \hat{n} \ = \ A_{\ell,m}^{\mathrm{a}}(r).
\Eeq
We could operate on $A_{\ell,m}^{\mathrm{a}}(r)$ in a number of ways, and reform a series similar to \eq{general-tensor} but using the appropriate bases.

\textsl{Grad} --- We can represent the gradient of a general tensor in similar form as \eq{general-tensor}
\Beq
\nabla \mathrm{T}(r,\theta,\phi) \ = \ \sum_{\mathrm{ab},\sigma\tau} \sum_{\ell,m}B_{\ell,m}^{\mathrm{ab}}(r) \,\mathcal{Q}_{\ell}(\sigma\tau,\mathrm{ab}) \, Y^{\overline{\sigma\tau}}_{\ell,m}(\theta,\phi)e(\sigma\tau)
\Eeq
where the radial coefficients, $B_{\ell,m}^{\mathrm{ab}}(r)$, are one rank higher than the input, $A_{\ell,m}^{\mathrm{b}}(r)$. That is,
\Beq
\label{tensor-grad}
B_{\ell,m}^{\mathrm{ab}}(r)  &=&  \xi^{\mathrm{\mathrm{a}}}_{\ell+\bar{\mathrm{b}}}  \, D_{\ell+ \bar{\mathrm{b}} }^{\mathrm{a}} \, A_{\ell,m}^{\mathrm{b}}(r),
\Eeq
and \eq{grad-unit} defines the $\xi^{\mathrm{a}}_{\ell+\bar{\mathrm{b}}}$ coefficients. 

This highlights two main points of this paper:
\begin{itemize}
\item
\textsl{\Eq{tensor-grad} implies that the gradient acts almost identically for vectors and scalars as for arbitrary-rank tensors.}
\item
\textsl{This is the main purpose for defining the $\mathcal{Q}_{\ell}$ transformations via \eq{Q-recursion}}
\end{itemize}

\textsl{Div} --- Computing the divergence requires a few relations between the $\mathcal{Q}_{\ell}$ transformations and $\xi$ factors. If $\sigma$, $\tau$ and $\mathrm{a}$, $\mathrm{b}$ are all single indices then contracting produces, 
\Beq
\label{ugly}
\sum_{\sigma,\,\tau}\delta_{\sigma,-\tau} \mathcal{Q}_{\ell}(\sigma \tau \kappa,\mathrm{abc})\ = \ \delta_{\mathrm{a},-\mathrm{b}} \sqrt{\frac{2 (\ell+\mathrm{b}+\bar{\mathrm{c}} )+1}{2 (\ell+ \bar{\mathrm{c}}  )+1}} \mathcal{Q}_{\ell}(\kappa,\mathrm{c}).
\Eeq
\Eq{ugly} is complicated to prove in general. It is straightforward to verify for any given specific rank. Also
\Beq
\label{xi-transformation}
\xi^{\mathrm{a} }_{\ell} = \xi^{-\mathrm{a}}_{\ell+\mathrm{a}} \sqrt{\frac{2(\ell+\mathrm{a})+1}{2\ell+1}},
\Eeq
which is a more straight forward fact than \eq{ugly}. Together, \eqs{ugly}{xi-transformation} imply 
\Beq
\label{contraction}
\sum_{\sigma,\,\tau}\delta_{\sigma,-\tau} \, \xi^{\mathrm{a}}_{\ell+\mathrm{b}+\bar{\mathrm{c}}} \mathcal{Q}_{\ell}(\sigma\tau\kappa, \mathrm{abc}) \ = \ \delta_{\mathrm{a},-\mathrm{b}} \, \xi^{\mathrm{b} }_{\ell+\bar{\mathrm{c}}} \mathcal{Q}_{\ell}(\kappa, \mathrm{c}).
\Eeq
\Eq{contraction} allows computing the divergences in terms of the radial derivative operators. 
Define  
\Beq
\label{another-general-tensor}
\mathrm{T}(r,\theta,\phi) &=& \sum_{\mathrm{ab},\sigma \tau} \sum_{\ell,m}A_{\ell,m}^{\mathrm{ab}}(r) \mathcal{Q}_{\ell}(\sigma\tau,\mathrm{ab}) Y^{\overline{\sigma\tau}}_{\ell,m}(\theta,\phi)e(\sigma\tau),
\\
\label{divergence-of-general-tensor}
\nabla \cdot \mathrm{T}(r,\theta,\phi) &=& \sum_{ \mathrm{b} , \tau} \sum_{\ell,m}B_{\ell,m}^{\mathrm{b}}(r) \mathcal{Q}_{\ell}(\tau, \mathrm{b}) Y_{\ell,m}^{\bar{\tau}}(\theta,\phi)e(\tau).
\Eeq
\Eq{another-general-tensor} is one rank higher than the output in \eq{divergence-of-general-tensor}; $\sigma$ and $\mathrm{a}$ represent the additional indices. The radially dependent components relate such that
\Beq
\label{general-div}
 B_{\ell,m}^{\mathrm{b}}(r)  \ = \ \sum_{\mathrm{a}} \xi^{\mathrm{a}}_{\ell+\bar{\mathrm{b}}}  \, D_{\ell+  \overline{\mathrm{ab}} }^{-\mathrm{a}} \, A_{\ell,m}^{\mathrm{ab}}(r),
\Eeq
where $-\mathrm{a}$ is the negative of index ``$\mathrm{a}$". Recall that ``$\mathrm{a}$'' represents a single index, while ``$\mathrm{b}$'' is possibly a multi-index. \Eq{general-div} is the tensor generalisation of \eq{div}.

\textsl{Laplacian} --- If we consider the action of $\nabla \cdot \nabla$ then the input and output components relate such that
\Beq
 B_{\ell,m}^{\mathrm{a}}(r)  \ = \  D_{\ell+  \bar{\mathrm{a}}  + 1}^{-}   \, D_{\ell+ \bar{\mathrm{a}} }^{+} \, A_{\ell,m}^{\mathrm{a}}(r) \ = \ D_{\ell+  \bar{\mathrm{a}}  - 1}^{+}   \, D_{\ell+ \bar{\mathrm{a}} }^{-} \, A_{\ell,m}^{\mathrm{a}}(r).
\Eeq
This shows the diagonalised Laplacian for general tensor fields. 

\section{General radial operations using Jacobi polynomials. }

At this point, our discussion of the geometric properties of three-dimensional space are complete enough to build a set of accurate and efficient numerical methods. The remainder of the paper follows much of V16. In the case of the 2D disk, it took much less background to get to the point where it is straightforward to compute numerical operations. In both V16 and the current paper the goal was to uncover a collection of spaces such as $\Reg{}$. These vector spaces are concrete enough to design ``bespoke'' computational algorithms.

\subsection{Radial basis and differentiation}

To show why the $D_{\ell}^{\pm}$ operators map between regularity spaces, we consider functions of the form 
\Beq
f_{\ell}(r) \ \equiv \  r^{\ell} \mathcal{F}(\z) \ \in \ \Reg{} , \quad \mathrm{where} \quad \z \ \equiv \ 2 r^{2} - 1,
\Eeq
and $\mathcal{F}(\z)$ is any differentiable function of $\z$. We make the specific convention of the even function $z$ of $r$ so that $ - 1 \le \z \le 1$. This is the same choice from V16 and is done to make calculations simpler for Jacobi polynomials. Therefore, consider the action of the differential operators on an element from $\Reg{}$,
\Beq
\label{D+(z)}
 \left[\frac{d}{dr} - \frac{\ell}{r}\right]r^{\ell} \mathcal{F}(\z) &=& 4 r^{\ell+1} \frac{d}{d\z} \mathcal{F}(\z) \\
 \label{D-(z)}
\left[\frac{d}{dr} + \frac{\ell+1}{r}\right]r^{\ell} \mathcal{F}(\z) &=& 2 r^{\ell-1} \left[ \left(\ell+\frac{1}{2}\right) + (1+\z)\frac{d}{d\z}\right]\mathcal{F}(\z) 
\Eeq
\Eqs{D+(z)}{D-(z)} are nothing more than the multiplication and chain rule of differentiation. From the leading-order behaviour on the right-hand side of \eqs{D+(z)}{D-(z)} we see that 
\Beq
D_{\ell}^{\pm} f_{\ell}(r)  \ \in \ \Reg{\pm 1},
\Eeq
which is a property shared with spherical Bessel functions.

Finding a good numerical method for the sphere is the same as finding a good basis for the regularity spaces. Jacobi polynomials $P$ are almost perfectly suited for the task. One reason is the following two differentiation formulae,
\Beq
\frac{d}{d\z} P_{n}^{(a,b)}(\z) &=& \frac{n+a+b+1}{2} P_{n-1}^{(a+1,b+1)}(\z) \label{Jacobi Dp} \\
\left[b + (1+\z)\frac{d}{d\z}\right]P_{n}^{(a,b)}(\z) &=& (n+b) P_{n}^{(a+1,b-1)}(\z), \label{Jacobi Dm}
\Eeq
for any parameters $a,b>-1$; see \cite{V16,DLMF}. We see that \eqs{D+(z)}{D-(z)} are consistent with \eqs{Jacobi Dp}{Jacobi Dm} if we choose the second Jacobi parameter 
\Beq
b \ \equiv \  \ell\ + \  \frac{1}{2}. \label{Jacobi b}
\Eeq
The first Jacobi parameter can take any value $a > -1$. \Eq{Jacobi b} contrasts with other choices in the literature; \eg \cite{livermore_jones_worland_2007} with $b = \ell- 1/2$.  The many properties of Jacobi polynomials allow many more useful relations than we will use here. The differential properties of the spin-weighted spherical harmonics in terms of $k_{\ell,s}^{\pm}$ are a direct consequence of \eq{Jacobi Dm} and a similar relation for $(a,b) \to (a-1,b+1)$. This is another main point of the paper:
\begin{itemize}
\item
\textsl{Carefully selecting Jacobi polynomial parameters provides big advantages in the efficiency of various differential operations.}
\end{itemize}

\Eq{Jacobi b} allows us to represent $f_{\ell}(r) \in \Reg{}$ as a series of Jacobi polynomials of type $(\alpha,\ell+ 1/2)$ and achieve simple differentiation formulae in the $r$ variable. We therefore define
\Beq
\label{Q poly def}
Q_{n}^{\alpha,\ell}(r) \ \equiv \    
\sqrt{\tfrac{2 \left(2n+\alpha+\ell+3/2\right) { {n+\alpha+\ell+1/2}\choose{\alpha} } }{ {{n+\alpha}\choose{\alpha}} }} \, 
 r^{\ell}  P_{n}^{(\alpha,\ell +1/2)}(2r^{2}-1).
\Eeq
We assume $\alpha > -1$; although it is possible to define negative parameter extensions. The parameter $n$ represents the degree of the Jacobi polynomials in terms of the $z$ argument. 
These polynomials are orthonormal on the weighed unit ball such that
\Beq
\int_{0}^{1} Q_{n}^{\alpha,\ell}(r) Q_{n'}^{\alpha,\ell}(r) \,(1-r^{2})^{\alpha}\, r^{2} \, \mathrm{d} r \ = \ \delta_{n,n'}.
\Eeq
Matsushima \& Marcus (1995) \cite{matsushima1995spectral} defined a parameterised basis set that is equivalent to \eq{Q poly def}. Their work found sparse differentiation formulae for these polynomials.  We find in V16 and in this current work that their matrices naturally factor into lower-bandwidth operators. We think the best term for the functions in \eq{Q poly def} are ``Generalised Spherical Zernike Polynomials''. This name reflects the fact that the $\alpha=0$  basis was discovered by Zernike for the disk in the 1930s \cite{zernike}. Also, the definition in \eq{Q poly def} is to the polynomials in V16 as spherical Bessel functions are to integer-index Bessel functions in 2D. 

After putting together the different operators, 
\Beq
\label{D+Q}
D^{+}_{\ell} Q_{n}^{\alpha,\ell}(r) &=&  Q_{n-1}^{\alpha+1,\ell+1}(r) \, \sqrt{2n 
   (2 n+2 \alpha+2 \ell+3)} \\
   \label{D-Q}
D^{-}_{\ell} Q_{n}^{\alpha,\ell}(r) &=&  Q_{n}^{\alpha+1,\ell-1}(r)  \, \sqrt{2
   (n+\alpha+1) (2 n+2 \ell+1)}.
\Eeq
The combination $2n+\ell$ represents the total  degree of the polynomial $Q_{n}^{\alpha,\ell}(r)$ in the $r$ variable. Therefore the right-hand side of both \eqs{D+Q}{D-Q} are polynomials of degree
\Beq
\deg\left[ D^{\pm}_{\ell} Q_{n}^{\alpha,\ell}(r) \right] \ = \ \deg\left[ Q_{n-1}^{\alpha,\ell+1}(r) \right] \ = \ \deg\left[ Q_{n}^{\alpha,\ell-1}(r) \right]  \ = \ \deg\left[Q_{n}^{\alpha,\ell}(r) \right] - 1. 
\Eeq

\Eqs{D+Q}{D-Q} are principal results of this paper. 
\begin{itemize}
\item
\textsl{The scalar Laplacian factors naturally into two first-order operators. These act as pure raising and lowering operators of the different Jacobi-polynomial parameters.}
\end{itemize}

\subsection{Conversion \& multiplication operations}

\Eqss{D+(z)}{D-Q} dictate how the second Jacobi index should behave. But just like in V16, we need to address how the first Jacobi index increments as we apply $D_{\ell}^{\pm}$; notice that $\alpha \to \alpha+1$ in both cases. The orthonormal basis in \eq{Q poly def} allows us to promote $\Reg{}$ to a collection of Hilbert spaces.  Define the spaces for each $\alpha>-1$, 
\Beq
\mathcal{H}^{\alpha}(\ell) \ \equiv \ \left\{ \, f \ \in \ \Reg{} \quad \mathrm{s.t.} \quad \int_{0}^{1} |f(r)|^{2} (1-r^{2})^{\alpha} r^{2} \dd{r} \ < \ \infty\,   \right\}.
\Eeq
From \eq{R +2 in R}, we still have 
\Beq
\mathcal{H}^{\alpha}(\ell+ 2) \ \subset \ \mathcal{H}^{\alpha}(\ell).
\Eeq
However, we also have the simple inequality, 
\Beq
(1 - r^{2})^{\alpha_{2}} \ \le \ (1 - r^{2})^{\alpha_{1}}, \quad \mathrm{for} \quad \alpha_{1} \ \le \ \alpha_{2},
\Eeq
which implies 
\Beq
\mathcal{H}^{\alpha}(\ell) \  \subset \  \mathcal{H}^{\alpha+1}(\ell).
\Eeq
We can now see the action of the differential operators more clearly
\Beq
D_{\ell}^{\pm} \ : \ \mathcal{H}^{\alpha}(\ell) \ \to \ \mathcal{H}^{\alpha+1}(\ell\pm 1).
\Eeq
\Eqs{D+Q}{D-Q} imply that this mapping is very sparse on the basis functions $Q_{n}^{\alpha,\ell}(r)$. 
Because successive spaces are embedded in each other, there is a well-behaved mapping between $\alpha \to \alpha+1$. Fortunately, Jacobi polynomials possess algebraic recursion formulae that allow incrementing the parameters and multiplying by $\z$-dependent factors \cite{V16,DLMF}. S\Beq
\label{Jacobi times}
\begin{array}{ccccc}
P_{n}^{(a,b)}(\z)             & = & \tfrac{n+a+b+1}{2n+a+b+1} P_{n}^{(a+1,b)}(\z)   & - & \tfrac{n+b}{2n+a+b+1} P_{n-1}^{(a+1,b)}(\z) \\ \\
\tfrac{1-\z}{2} P_{n}^{(a,b)}(\z) & = &   - \tfrac{n+1}{2n+a+b+1} P_{n+1}^{(a-1,b)}(\z) & + & \tfrac{n+a}{2n+a+b+1} P_{n}^{(a-1,b)}(\z)\\ \\
 P_{n}^{(a,b)}(\z) & = & \tfrac{n+a+b+1}{2n+a+b+1} P_{n}^{(a,b+1)}(\z)  & + & \tfrac{n+a}{2n+a+b+1} P_{n-1}^{(a,b+1)}(\z)  \\ \\ 
\tfrac{1+\z}{2} P_{n}^{(a,b)}(\z) &=&   \tfrac{n+1}{2n+a+b+1} P_{n+1}^{(a,b-1)}(\z) &+& \tfrac{n+b}{2n+a+b+1} P_{n}^{(a,b-1)}(\z).
\end{array}
\Eeq

After including normalisation factors, we find the Jacobi recursion operators become various 
multiplication operators acting on the $Q_{n}^{\alpha,\ell}(r)$ functions\footnote{\Eq{Jacobi times} also implies the $\sin$/$\cos$ multiplication operators acting on spin-weighted spherical harmonics in \eqss{cos times}{sin times}; after taking into account the different normalisation.},  
\Beq
\label{multiplication on Q}
\begin{array}{ccccc}
Q_{n}^{\alpha,\ell}(r)                & = &  Q_{n}^{\alpha+1,\ell}(r)\, \mathrm{A}_{n}  & + & Q_{n-1}^{\alpha+1,\ell}(r)\, \mathrm{B}_{n-1} \\ \\
(1-r^2) Q_{n}^{\alpha+1,\ell}(r) & = &   Q_{n+1}^{\alpha,\ell}(r)\, \mathrm{B}_{n} &+& Q_{n}^{\alpha,\ell}(r)\, \mathrm{A}_{n}  \\ \\
 r\,Q_{n}^{\alpha,\ell}(r)  & = & Q_{n}^{\alpha,\ell+1}(r)\, \mathrm{C}_{n}  &+& Q_{n-1}^{\alpha,\ell+1}(r)\, \mathrm{D}_{n-1}  \\ \\ 
 r\,Q_{n}^{\alpha,\ell+1}(r) &=&   Q_{n+1}^{\alpha,\ell}(r)\, \mathrm{D}_{n} &+& Q_{n}^{\alpha,\ell}(r)\,  \mathrm{C}_{n}
\end{array}
\Eeq
with the coefficients defined by, 
\Beq
\begin{array}{ccccccc}
\mathrm{A}_{n} &=& \sqrt{\!\tfrac{(n+\alpha+1) (n+\alpha+\ell +3/2)}{(2 n+\alpha+\ell +3/2) (2n+\alpha+ \ell +5/2)}}, & & \mathrm{B}_{n} &=& - \,\sqrt{\!\tfrac{ (n+1) (n+\ell +3/2)}{( 2 n+\alpha+\ell +5/2) (2 n+\alpha+\ell
   +7/2)}}\\ \\
   \mathrm{C}_{n} &=& \sqrt{\!\tfrac{\left(n+\ell +3/2\right)\left(n+\alpha+\ell +3/2\right)}{
\left(2 n+ \alpha+\ell +3/2\right)\left(2 n+\alpha+\ell+5/2\right)}}, & & \mathrm{D}_{n} &=& \sqrt{\!\tfrac{(n+1) (n+\alpha+1)}{\left(2n+\alpha+\ell +5/2\right)
   \left(2 n+\alpha+\ell +7/2\right)}} \, .
\end{array}\label{ABCD defs}
\Eeq
The $\alpha$ and $\ell$ dependence in the coefficients in \eq{ABCD defs} is the same for each element; only $n$ changes. In all cases $Q_{-1}^{\alpha,\ell}(r) = 0$ by definition.

We expand functions in $\mathcal{H}^{\alpha}(\ell)$ in terms of the orthonormal basis, 
\Beq
f(r) \ = \ \sum_{n\ge0} Q_{n}^{\alpha,\ell}(r)\, F_{n}. \label{spectral sum}
\Eeq
We compute the spectral coefficients via projection 
\Beq
F_{n} \ = \ \int_{0}^{1} Q_{n}^{\alpha,\ell}(r) \,f(r)\, (1-r^{2})^{\alpha}\, r^{2} \dd{r}. \label{spectral projection}
\Eeq
The coefficients, $F_{n}$, depend implicitly on $\alpha$ and $\ell$, but we omit the labels for simplicity.
For a finite number of modes in the sum in \eq{spectral sum}, we can compute \eq{spectral projection} using Gaussian quadrature for Jacobi polynomials. We want a simple way to implement various operators on spectral coefficient data. Therefore, we define a \textit{row} vector of basis elements, 
\Beq
Q^{\alpha,\ell}(r) \ = \ \left[ \,Q_{0}^{\alpha,\ell}(r) ,\,  Q_{1}^{\alpha,\ell}(r), \,Q_{2}^{\alpha,\ell}(r),\, \ldots \, \right].
\Eeq
We use a tilde to denote the following \textit{column} vector, \ie
\Beq
\widetilde{Q}^{\alpha,\ell}(r) \ \equiv \  r^{2} (1-r^{2})^{\alpha}\, Q^{\alpha,\ell}(r)^{\top}. \label{dual Q}
\Eeq
\Eq{dual Q} gives the dual basis in the sense that, 
\Beq
Q^{\alpha,\ell}(r) \cdot \widetilde{Q}^{\alpha,\ell}(r')  \ = \delta(r-r'), \quad \mathrm{and} \quad \int_{0}^{1}  \widetilde{Q}^{\alpha,\ell}(r)\, Q^{\alpha,\ell}(r) \dd{r} \ = \ I,
\Eeq
where $\delta(r-r')$ is the Dirac-$\delta$ distribution, and $I$ is the identity matrix acting on the $n$ index. These definitions imply that the sum and projection in \eqs{spectral sum}{spectral projection} become,
\Beq
f(r) \ = \ Q^{\alpha,\ell}(r) \cdot {F}  \quad \mathrm{and} \quad F = \int_{0}^{1} \widetilde{Q}^{\alpha,\ell}(r) \,f(r) \dd{r} . \label{projection}
\Eeq
We make the choice of the row/column convection so our state vector of spectral coefficients, $F$, is represented as a column vector.  

With the above definitions, the algebraic and differential operators act on the basis elements in the following way,
\Beq
\begin{array}{ccccccc} 
Q^{\alpha,\ell}(r) &=& Q^{\alpha+1,\ell}(r)\, C_{\alpha,\ell}, & {} &  (1-r^2)\, Q^{\alpha+1,\ell}(r) & = &   Q^{\alpha,\ell}(r)\, {C}^{\top}_{\alpha,\ell} 
\\
& & & {} & & & 

\\
r\,Q^{\alpha,\ell}(r) &=& Q^{\alpha,\ell+1}(r)\, R_{\alpha,\ell}, & {} & r \, Q^{\alpha,\ell+1}(r) &=&  Q^{\alpha,\ell}(r)\, {R}^{\top}_{\alpha,\ell}
\\ 
& & & {} & & &
\\
\left[\frac{d}{dr} - \frac{\ell}{r} \right] Q^{\alpha,\ell}(r) &=& Q^{\alpha+1,\ell+1}(r) D^{+}_{\alpha,\ell}, & {} & \left[\frac{d}{dr} + \frac{\ell+1}{r} \right] Q^{\alpha,\ell}(r) &=& Q^{\alpha+1,\ell-1}(r) D^{-}_{\alpha,\ell}\,.
\end{array}\label{abstract operators}
\Eeq
That is, starting with $f \in \mathcal{H}^{\alpha}(\ell)$:
\Beq
\begin{array}{ccccccc}
f(r) \in \mathcal{H}^{\alpha+1}(\ell) &  \longleftrightarrow & C_{\alpha,\ell} \cdot F, &  & (1-r^2) f(r) \in \mathcal{H}^{\alpha-1}(\ell) &  \longleftrightarrow & C^{\top}_{\alpha-1,\ell} \cdot F \\ \\ 
r\, f(r) \in \mathcal{H}^{\alpha}(\ell+1) &  \longleftrightarrow & R_{\alpha,\ell} \cdot F, &  & r\, f(r) \in \mathcal{H}^{\alpha}(\ell-1) &  \longleftrightarrow & R^{\top}_{\alpha,\ell-1} \cdot F \\ \\
D_{\ell}^{+}f(r) \in \mathcal{H}^{\alpha+1}(\ell+1) & \longleftrightarrow & D^{+}_{\alpha,\ell} \cdot F, &  & D^{-}_{\ell}f(r) \in \mathcal{H}^{\alpha+1}(\ell-1) & \longleftrightarrow & D^{-}_{\alpha,\ell} \cdot F\, ,
\end{array}
\Eeq
where in the last set of correspondences, the $D_{\ell}^{\pm}$ operators acting on function represents the differential operators defined in \eq{D+,D-}, and $D_{\alpha,\ell}^{\pm}$ represent sparse matrices defined via \eqs{D+Q}{D-Q}. 

In \eq{abstract operators}, $C_{\alpha,\ell}, R_{\alpha,\ell}$, and $D^{\pm}_{\alpha,\ell}$ all represent sparse matrices that act from the left on spectral coefficients. The matrix $D_{\alpha,\ell}^{-}$ is diagonal, $D_{\alpha,\ell}^{+}$ has a single band immediately above the diagonal, and each $C_{\alpha,\ell}$ and $R_{\alpha,\ell}$ matrix has entries along the diagonal and one the band immediately above. The following duality is an interesting pattern of the  various operators: 
\begin{enumerate}
\item
Operators that are local in physical space, \ie $1$, $r$, $1-r^{2}$, couple two neighbouring modes in the $n$ index. 
\item
Operators coupling nearby points in physical space, \ie $D_{\ell}^{\pm}$, do not couple modes in the $n$ index. 
\end{enumerate}
Operators that are local in physical space are non-local in coefficient $n$ space, and the converse is also true.

Finally, for computing the solution to boundary-value problems we need the restriction operator acting at $r=1$. This is as simple as,
\Beq
f(r=1) \ = \ Q^{\alpha,\ell}(r=1) \cdot F 
\Eeq
For Jacobi polynomials, 
\Beq
P_{n}^{(a,b)}(\z=1) \ = \ {{n+a}\choose{a}} \ \equiv \ \frac{\Gamma(n+a+1)}{\Gamma(a+1)\, n!}.
\Eeq
After adjusting for the normalisation,
\Beq
\label{Q-r=1}
Q_{n}^{\alpha,\ell}(r=1) \ \equiv \    
\sqrt{2 \left(2 n+ \alpha+\ell+3/2\right) { {n+\alpha+\ell+1/2}\choose{\alpha} } {{n+\alpha}\choose{\alpha}}},
\Eeq
For $n\gg1$, 
\Beq
Q_{n}^{\alpha,\ell}(r=1) \ \sim \  \frac{2\, n^{\alpha+1/2}}{\Gamma(\alpha+1)}  , \quad \mathrm{as} \quad n \ \to \ \infty \label{large-n restriction},
\Eeq
which gives an idea of how well a function in $\mathcal{H}^{\alpha}(\ell)$ must behave to have a meaningful restriction to the boundary. \Eq{large-n restriction} implies that restriction is likely to make sense for the vast majority of circumstances encountered. 

As with V16, all the operators except $C_{\alpha,\ell}^{\top}$ are essential for constructing algebraic systems out of general linear PDEs. However, since $1-r^{2} = 0$ at $r=1$,  the $C_{\alpha,\ell}^{\top}$ matrices are  useful for creating basis sets that automatically satisfy Dirichlet boundary conditions. 

From these mappings, we can easily spot a few useful relationships. For example, 
\Beq
(1-r^{2}) \, f(r) \in \mathcal{H}^{\alpha}(\ell) &  \longleftrightarrow & C^{\top}_{\alpha,\ell}\,C_{\alpha,\ell}\cdot F \ = \ C_{\alpha-1,\ell}\,C_{\alpha-1,\ell}^{\top}\cdot F,
\Eeq
and
\Beq
r^{2} \, f(r) \in \mathcal{H}^{\alpha}(\ell) &  \longleftrightarrow & R^{\top}_{\alpha,\ell}\,R_{\alpha,\ell}\cdot F \ = \ R_{\alpha,\ell-1}\,R_{\alpha,\ell-1}^{\top}\cdot F.
\Eeq
All of these are manifestly positive definite both because $0 \le 1-r^2 \le 1$, and because $M^{\top}M$ is self-adjoint for any real matrix. Furthermore (for example),  
\Beq
C^{\top}_{\alpha,\ell}\,C_{\alpha,\ell} + R^{\top}_{\alpha,\ell}\,R_{\alpha,\ell} \ = \  I.
\Eeq

\section{Discussion}

At this point we have derived methods for solving a large number of possible linear PDEs in spherical coordinates. The principles behind our methods are two-fold:
\begin{enumerate}
\item First we leverage the geometry as much as possible to decompose a high-dimensional problem into a (possibly large) collection of decoupled one-dimensional problems. This stage relies on properties specific to the sphere and ball in three dimension. The methods are exact, and could be considered analytical as opposed to numerical or computational. The computational effectiveness of the exact methods relies on their expressiveness; which is to say the amount of useful information contained in an easily computable expression. We believe that the formulae provided in \S1--4 meet this criteria. In the Part-II companion paper, we discuss in detail how we create the operators and data structures in easy-to-use code, and apply those operators to interesting PDE problems in the sphere.

\item After distilling the geometric content from three-dimensional polar-coordinate tensors, we then focus on the numerical task of solving equations within the regularity classes $\Reg{}$. This is a somewhat seperate question from how we arrive at the classes in the first place. Functions with a particular behaviour near a coordinate singularity can be represented either with a specially designed basis, or with something more general purpose. In principle, virtually any numerical method employed in the solution of PDEs can solve equations on $\Reg{}$. Some of these methods do much better than others. In \S4 we provide a class of methods (with a single free parameter) that conform specifically to the type of differential operators that commonly act in the sphere. These operators are not exact eigenfunctions of the given operators. In this sense they make a general-purpose method. On the other hand, explicitly considering geometry allows for a set of maximally sparse matrices that compute common operators accurately and fast on general data. Our Jacobi-polynomial basis is physics agnostic, but adheres to the geometry. We could however, choose an even more general-purpose methods. We make some comments on this below.
\end{enumerate}

\subsection{Matrix-system construction}

There is a large literature discussing the advantages and disadvantages of different methods for solving equations of the form (for example) of
\Beq
D_{\ell+1}^{-} D_{\ell}^{+} f \ = \ - \kappa^{2} f , \quad f(r=1) = 0. \label{Laplace eigenvalues}
\Eeq
The exact solutions to \eq{Laplace eigenvalues} are spherical Bessel functions
\Beq
f_{n}(r) \ \propto \ j_{\ell}(\kappa_{n} r) \, \quad \mathrm{with}\quad j_{\ell}(\kappa_{n}) \ = \ 0.
\Eeq
As with trigonometric functions on the interval, Bessel functions make excellent test problems for orthogonal polynomials. The reason is that like $\sin/\cos$, Bessel functions are very much not finite polynomials. Highly accurate solutions are not a \textit{fait accompli}. A method that does well on this kind of simple test problem is likely a good general-purpose method. 
The companion paper (Part-II) discusses the details of much more complicated numerical examples. For now, we only want to show a simple example for the purpose of demonstrating how to construct systems of linear equations from our method. This example also gives some context for how the same systems are constructed and solved in competing general-purpose methods.

The key question is: how should we turn \eq{Laplace eigenvalues} into a set of matrices that we can solve numerically? 

First we pick a value for $\alpha > -1  $ for our primary basis set. The values $\alpha=0$ and $\alpha = -1/2$ are both important special cases.  The case $\alpha=0$ corresponds to the geometric weight in 3D, $r^{2} \dd{r}$. The case $\alpha = -1/2$ corresponds to $r^{2} \dd{r}/\sqrt{1-r^{2}}$. This second case has the disadvantage of not using the proper geometric measure. It is possible to circumvent this minor problem by making small adjustments to the quadrature weights when computing integration over the spherical volume.  In any case,
\Beq
f = {Q}^{\alpha,\ell}\cdot F.
\Eeq
In the terminology of finite-element methods, we choose $Q^{\alpha,\ell}$ as a basis of \textit{trial functions}. 
After applying the left-hand side of \eq{Laplace eigenvalues}, we find
\Beq
{Q}^{\alpha+2,\ell}\cdot D_{\alpha+1,\ell+1}^{-}\,D_{\alpha,\ell}^{+} F \ = \ -\kappa^{2}\, {Q}^{\alpha,\ell}\cdot F. \label{LEP step-1}
\Eeq
The problem with \eq{LEP step-1} is that the left-hand side is now in the $\alpha+2$ basis, while the right-hand side is still in the $\alpha$ basis. But we can \textit{convert} the basis on the right up two levels.
\Beq
{Q}^{\alpha+2,\ell}\cdot \left [ D_{\alpha+1,\ell+1}^{-}\,D_{\alpha,\ell}^{+} \ + \ \kappa^{2}\,  C_{\alpha+1,\ell}\,C_{\alpha,\ell}  \right]  F \ = \ 0.\label{LEP step-2}
\Eeq
We can now project \eq{LEP step-2} on the the $Q^{\alpha+2,\ell}$ basis. In the language of finite-element methods, this is the basis of \textit{test functions}. In general, the trial functions and test functions can lie in seperate spaces. Many spectral methods produce dense matrices because the trial and test spaces are considered to be the same. This may seem like a simplification, but simplifying the basis set at the cost of sparseness in the matrices can lead to a large amount of computational cost. 

We also must include boundary conditions,
\Beq
{Q}^{\alpha,\ell}(r=1)\cdot F \ = \ 0. \label{LEP step-3}
\Eeq
Currently, \eq{LEP step-2} is \textit{upper-triangular} in terms of the array of variable, $F$. Upper-triangularity implies that the system is closed if we truncate with $N$ equations for $N$ variables \textit{without considering the boundary conditions}. The boundary condition \eq{LEP step-3} makes for $N+1$ equations for $N$ variables, and the system is over constrained. We could solve the problem by truncating \eq{LEP step-2} with only $N-1$ equations for $N$ variables. But there is a more general method that can accomplish the same thing as well.

We employ a $\tau$-method correction \cite{Lanczos38, Ortiz69, boyd_book} in \eq{LEP step-2}. That is, we include a right-hand side correction of the form 
\Beq
{Q}^{\alpha+2,\ell}\cdot \left [ D_{\alpha+1,\ell+1}^{-}\,D_{\alpha,\ell}^{+} \ + \ \kappa^{2}\,  C_{\alpha+1,\ell}\,C_{\alpha,\ell}  \right]  F \ = \ - \tau \, P^{\ell}_{N} \label{LEP step-4}
\Eeq
In \eq{LEP step-4} $P^{\ell}_{N}(r) \in \Reg{}$ is any \textit{specific} polynomial with
\Beq
\mathrm{deg} \left[ P^{\ell}_{N}(r) \right] \ = \  \ell+ 2 N.
\Eeq
The variable $\tau$ is a single additional degree of freedom that allows imposing the boundary conditions. The choice of specific $P^{\ell}_{N}(r)$ determines the specific truncation of the system.  If we choose, $P^{\ell}_{N}(r) = Q_{N}^{\alpha+2,\ell}(r)$, then the $\tau$-correction is equivalent to dropping the last equation from the system a replacing it with the boundary conditions.  But we can make more general choices that can affect the accuracy of the solution. We explore the consequences of different $\tau$-correction choices in Part-II.

We can now project \eq{LEP step-3} onto the ${Q}^{\alpha+2,\ell}$ basis. The entire system written in block form is
\Beq
\left [  
\begin{array}{  c c c  }
\hspace{0.1cm}
D_{\alpha+1,\ell+1}^{-}\,D_{\alpha,\ell}^{+} \ + \ \kappa^{2}\,  C_{\alpha+1,\ell}\,C_{\alpha,\ell} & & \Pi_{\,N}^{\alpha,\ell} \hspace{0.1cm} \vspace{0.2cm} \\ 
 \hspace{0.1cm} {Q}^{\alpha,\ell}(r=1) & & 0 \hspace{0.1cm} 
\end{array} 
 \right] . \left[ \begin{array}{c} F  \vspace{0.2cm} \\ \tau \end{array} \right] \ = \ \left[ \begin{array}{c} 0  \vspace{0.2cm} \\ 0 \end{array} \right] . \label{LEP step-5} 
\Eeq
where 
\Beq
\Pi_{\,N}^{\alpha,\ell} \ \equiv \ \int_{0}^{1}  \widetilde{Q}^{\alpha+2,\ell}(r)\, P^{\ell}_{N}(r) \dd{r}
\Eeq
is the column vector of the $\tau$-correction polynomial projected against the test basis.

If $P^{\ell}_{N}(r) \ = \ Q_{N}^{\alpha+2,\ell}(r)$, then 
\Beq
\Pi_{\,N}^{\alpha,\ell} \ = \ E_{N} \ \equiv \ \left[ 0, \ldots ,\, 0, \, 1 \right]^{\top} 
\Eeq
 is a column of all 0s, except for a single 1 in the last entry. However, we could choose $P^{\ell}_{N}(r) \ = \ Q_{N}^{\alpha,\ell}(r)$, in which case, 
\Beq
\Pi_{\,N}^{\alpha,\ell} \ = \   C_{\alpha+1,\ell}\,C_{\alpha,\ell}\, E_{N}.
\Eeq
In this case, $\Pi_{\,N}^{\alpha,\ell}$ only has three entries (rather than one) and the system is still sparse. We can pick $P_{N}^{\ell}(r)$ to maximise accuracy (under some constraints), and then we can pick $\alpha$ to maximise the sparsity (under some other constraints). It is also possible with a small number of matrix operations to reduce the $\tau$ column to a single entry and therefore be able to decouple it entirely form the rest of the system. 

\subsection{Automatic boundary conditions}

There is also another quite elegant method for implementing the boundary conditions, which we turn to now.
We note that the adjoint of the conversion matrix $C$ has the effect of multiplying by $(1-r^2)$. We can therefore propose that 
\Beq
F \ = \ C^{\top}_{\alpha,\ell}\, G,\label{basis-recombo}
\Eeq
where $G$ is now an unknown column vector. \Eq{basis-recombo} implies that 
\Beq
f(r) \ = \ {Q}^{\alpha,\ell}(r) \cdot F \ = \ {Q}^{\alpha,\ell}(r) \cdot C^{\top}_{\alpha,\ell}\, G \ = \ (1-r^{2})\, {Q}^{\alpha+1,\ell}(r) \cdot G  \label{f(1) like zero}
\Eeq
This implies that $ f(r=1) =  0$ automatically. \Eqs{basis-recombo}{f(1) like zero} imply that it is sometimes advantageous to use a non-orthogonal basis for trial functions.  We now only need to apply the appropriate matrix on the right of the system and we are done:
\Beq
\left[ D_{\alpha+1,\ell+1}^{-}\,D_{\alpha,\ell}^{+}  + \kappa^2 \, C_{\alpha+1,\ell}\,C_{\alpha,\ell}  \right] C^{\top}_{\alpha,\ell}\, G \ = \ 0. \label{LEP step-6}
\Eeq
Livermore (2010) \cite{Livermore} pointed out some notable algebraic properties of Jacobi polynomials. In particular their ability to easily satisfy Galerkin boundary conditions and still maintain some kind of orthogonality. We believe that \eq{f(1) like zero} lies at the root of that ability. 
We could make small adjustments to the matrix in \eq{basis-recombo} if we want a different boundary condition other than Dirichlet (see \S5.2 of V16, and also \cite{doha2006efficient,julien_watson_2009,Livermore,li2010optimal,Livermore_Ierley_2010}). The matrices in \eq{LEP step-6} have one band on the lower diagonal, and two bands above; there are no dense rows. This approach is called the ``Dirichlet basis'' method, the ``Galerkin basis'' method, or ``basis recombination''; depending on the context.  

With a few more moving parts, it is possible to combine the $\tau$-correction method with the automatic boundary conditions. Even if the $\tau$ column is dense, it is possible to row-precondition the system in a similar way as for the automatic boundary conditions, and increase the matrix bandwidth by one more band.  This leads us to highlight another important theme of our work:
\begin{itemize}
\item
\textsl{The \textit{sparsity} and the \textit{accuracy} of a given system are two seperate issues that are all-too-often intertwined. We can make a system both minimally banded and also choose an arbitrary (polynomial) space to absorb the truncation error.}
\item
\textsl{Sparsity is advantageous for the speed of computing solutions, numerical well-conditioning, and ease of automatic system construction.}
\item
\textsl{Slightly different $\tau$-correction terms can produce large differences in the overall accuracy of a scheme, independent of the sparsity.}
\item
\textsl{A specific choice of $\tau$-correction calls into existence an \textit{exact} finite solution in terms of polynomials. At this point, \textit{any} basis of polynomials will solve to problem; some of them with many more numerical operations required than others.} 
\end{itemize}

\subsection{Comparison to other work}

It is also helpful to compare the basis sets and operators described in \S4 to some other tools that exist for problems in the full sphere. It is also useful to compare our methods to similar schemes using Chebyshev polynomials in domains without coordinate singularities (\eg Cartesian, and cylindrical annulus).

Other than the generality in the choice of $\tau$-correction, \eqs{LEP step-3}{LEP step-4} follow the same approach that we took in V16.  A number of other authors have taken similar approach in one-dimension with Chebyshev polynomials \cite{zebib_1984,coutsias_hagstrom_1996,boyd_book,greengard_1991,doha2002efficient,julien_watson_2009,doha2002efficient,muite_2010,olver2013fast,viswanath_2015,trefethenspectralmethods}. \Eq{LEP step-3} is the closest in methodology to the ``ultra-spherical method" of Olver and Townsend (2012) \cite{olver2013fast}. Our method also shares some similarities to the ``quasi-inverse'' method of Julien and Watson (2009) \cite{julien_watson_2009}.  In the latter case, the ``quasi-inverse'' implied finding a pre-conditioning operator that renders the equations sparse. The pre-conditioning operator followed systematically from a three-term recursion formulae for the derivative of any basis set.  It turns out that Jacobi polynomials are the unique class of basis that satisfies the properties require for the quasi-inverse. We also know now that for Jacobi polynomials and PDEs, the preconditioning matrices will always relate to a basis-conversion operator in some way.

We can truncate \eq{LEP step-3} in many different ways without using a $\tau$ correction. However, all methods of truncation can be mapped into a $\tau$-method with some kind of equation correction; see 
\cite{Daou_Ortiz}\footnote{The ``$\infty$-dimensional'' method of Olver and Townsend (2014) \cite{olver2014practical} is the only scheme that differs much from this approach. In the $\infty-$dimensional case, the idea is to preform sparse orthogonal matrix operations to render part of the system upper-triangular before deciding the degree of truncation based on an adaptive error estimate. The method is infinite in the sense that it uses a greedy algorithm in the forward orthogonal transformation stage.}. For many examples, truncating the system and replacing a row with \eq{LEP step-4} works quite well. However, it is not automatic that all truncations work equally well \cite{Charalambides_Waleffe_2008,Dawkins_etal_1998,gardner_etal_1989,mcfadden_etal_1990}. A full discussion of which method works best is an important topic that is outside the scope of our work here.

\subsubsection{Jones-Worland Polynomials} 

In recent years, a number of author have used a basis set that they have termed Jones-Worland polynomials \cite{livermore_jones_worland_2007,Hollerbach+2013,Marti_Jackson,Marti-benchmark}. This basis relates in a clear way to the basis sets that we use. The individual basis elements $W$ are defined in terms of Jacobi polynomials $P$ (up to normalisation) such that 
\Beq
W_{n}^{\ell}(r) \ \propto \ r^{\ell}\, P_{n}^{(-\tfrac{1}{2},\, \ell-\tfrac{1}{2})}(2r^{2}-1) \label{Worland}
\Eeq
This basis is similar to our $\alpha=-1/2$ basis. But note the important difference that the second Jacobi parameter is $\ell-1/2$, rather than $\ell+1/2$ from \eq{Jacobi b}, which allows for the basis to most efficiently conform to spherical geometry. The proponents of \eq{Worland} say that the basis set has the advantage of being the closest analogue of Chebyshev polynomials with as close-to-uniform oscillation as possible. 

In some of our numerical tests in Part-II, we find some advantage from using $\alpha = -1/2$ over $\alpha=0$; but not in every situation. When $\alpha = -1/2$ is better, we find that almost all of the gains come from controlling behaviour near the outer boundary. Half-integer differences between values of Jacobi parameters are not accessible with sparse transformations; integer differences allow sparse conversion.\footnote{However, it is possible to make the conversion fast in a highly accurate, albeit, approximate way. This technique is useful for computing transforms between a spatial grid and spectral coefficients \cite{{slevinsky_2018,Slevinsky_SHT}}.}  We see from \eq{large-n restriction} that $\alpha=-1/2$ does produce uniformity in the boundary values of the basis functions. However integer differences for the geometrically motivated parameter do not make much difference in the numerical properties of solutions.

There is a simple conversion between the polynomials defined in Livermore, Jones \& Worland (2007)  \cite{livermore_jones_worland_2007} and our definition of a Generalised Spherical Zernike basis set. For $\alpha =-1/2$, 
\Beq
W^{\ell}(r) \ = \ r\, Q^{\alpha,\ell-1}(r) \ = \ Q^{\alpha,\ell}(r)\, R_{\alpha,\ell-1}.
\Eeq
The best this basis can do in terms of sparsity would be to replace our operators in \eq{LEP step-3} such that
\Beq
D_{\alpha+1,\ell+1}^{-}\,D_{\alpha,\ell}^{+} + \kappa^2 \, C_{\alpha+1,\ell}\,C_{\alpha,\ell}  \ \to \ \left [D_{\alpha+1,\ell+1}^{-}\,D_{\alpha,\ell}^{+} + \kappa^2 \, C_{\alpha+1,\ell}\,C_{\alpha,\ell}\right] R_{\alpha,\ell-1}.\label{worland-problem}
\Eeq
Previous work using this basis have always used a dense formulation such that, 
\Beq
D_{\alpha+1,\ell+1}^{-}\,D_{\alpha,\ell}^{+} + \kappa^2 \, C_{\alpha+1,\ell}\,C_{\alpha,\ell}  \ \to \ R_{\alpha,\ell-1}^{-1} \, C_{\alpha,\ell}^{-1}\, C_{\alpha+1,\ell}^{-1}\, D_{\alpha+1,\ell+1}^{-}\,D_{\alpha,\ell}^{+} \, R_{\alpha,\ell-1} + \kappa^2.
\Eeq
The inverse operators make the problem dense. Even in the sparse \eq{worland-problem}, there is no need to add an additional matrix band. If the Jones-Worland basis gives different numerical results for a specific problem, it is always possible to make small alterations to the last entries of a truncated matrix to render the two methods exactly equivalent. This could be done with a small change in a $\tau$-correction polynomial; \ie
\Beq
\Pi_{\,N}^{\alpha,\ell} \ = \   C_{\alpha+1,\ell}\,C_{\alpha,\ell}\,R_{\alpha,\ell-1} \, E_{N}, \quad \mathrm{with} \quad \alpha = - 1/2.
\Eeq
This column contains four non-zero entries. 

Figure 1 shows the relative error of the first solution to the spherical Bessel function problem using three different methods. The green and red lines use respective $\tau$ columns 
\Beq
\Pi_{\,N}^{0,\ell} \ = \   C_{1,\ell}\,C_{0,\ell}\, E_{N}, \quad \mathrm{or} \quad \Pi_{\,N}^{-1/2,\ell} \ = \   C_{3/2,\ell}\,C_{-1/2,\ell}\, R_{-1/2,\ell-1} E_{N}.\label{tau-columns}
\Eeq
The plots compute the relative error from solving for the first root of a spherical Bessel function, 
\Beq
j_{\ell}(\kappa_{0,\ell}) \ = \ 0, \quad \mathrm{where} \quad \kappa_{0,\ell} \approx \begin{cases} 5.7634591968945498, & \ell= 2 \\
72.199780933233055, &    \ell= 64.
\end{cases}
 \Eeq
We can see that the different methods produce different slopes of the error curves. Small changes can produce large difference at low-to-moderate resolutions. We find this is very important in Part-II when studying the some of the finer details of magnetic dynamos. Both methods from \eq{tau-columns} do much better than Chebyshev methods at the same number of degrees of freedom.

\begin{figure}
  \centerline{\includegraphics[width=16cm]{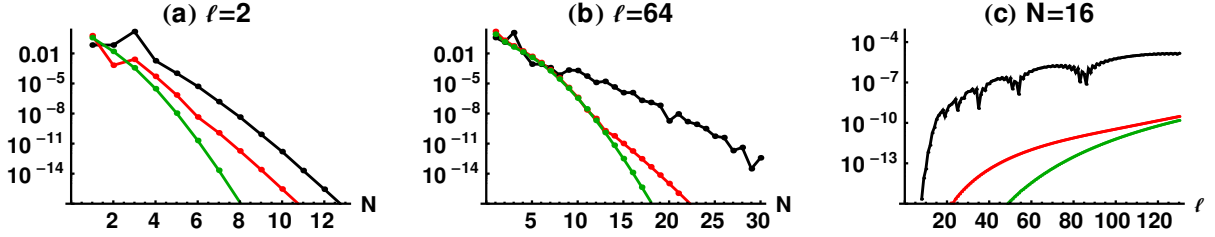}}
  \caption{Errors from different numerical methods: The plots all show the relative error in the fundamental wavenumer for the spherical Bessel function problem; $|\kappa_{0,\ell}^{N} - \kappa_{0,\ell}|/\kappa_{0,\ell}$. Plots (a) and (b) show the error at fixed spherical harmonic degree, $\ell$ as a function of the total degrees of freedom, $N$. Plot (c) shows the error at a fixed number of degrees of freedom as a function of regularity,. $\ell$. In each plot, the green line (lowest error) shows result from a $\tau$-correction with $\alpha=0$. The red line (medium error) shows the error from a ``Jones-Worland'' $\tau$-correction. The black line (largest error) shows the error from a sparse Chebyshev solution. }
\label{fig:eigenvalue errors}
\end{figure}

\subsubsection{The recursion formulae of Matsushima \& Marcus}

Matsushima \& Marcus (1995) \cite{matsushima1995spectral} found a basis with sparse differentiation operators. Their basis sets were equivalent to our basis set for specific parameters. We show here how their operators relate.

Matsushima \& Marcus used operators on the form
\Beq
 r \frac{d}{dr} \ = \ r D_{\ell}^{+} + \ell \ = \ r D_{\ell}^{-} -  (\ell+1) \ = \ \frac{1}{2}\left( r D_{\ell}^{+} + r D_{\ell}^{-} -  1 \right).
\Eeq
It is possible to construct sparse matrix systems by multiplying through by an appropriate power of $r$. For example,   
\Beq
r^{2}\, D_{\ell+1}^{-} D_{\ell}^{+} \ = \ \vartheta ( \vartheta + 1 ) - \ell( \ell+ 1), \quad \mathrm{where} \quad \vartheta \ \equiv \ r \frac{d}{dr}. 
\Eeq
Acting on a basis set,
\Beq
2 r \frac{d}{dr} Q^{\alpha,\ell}(r)\, C_{\alpha-1,\ell}  \ = \  Q^{\alpha,\ell}(r) \left( R_{\alpha,\ell-1} D_{\alpha-1,\ell}^{-} + R_{\alpha,\ell}^{\top}\, D_{\alpha-1,\ell}^{+} - C_{\alpha-1,\ell} \right)
\Eeq
That is, 
\Beq
r \frac{d}{dr}  \ : \ \mathcal{H}^{\alpha}(\ell) \ \to \ \mathcal{H}^{\alpha+1}(\ell).
\Eeq
Matsushima \& Marcus also found other related operations, \eg $(1-r^{2})$, but not $r$ individually. The operations are presented in terms of explicit recursion formulae for the polynomials \textit{within} a regularity space. The most important distinction was that there are not operators that can map \textit{between} regularity spaces.  This not only makes systems more difficult to construct automatically, but also significantly more dense.  

\subsection{Chebyshev solution}

Several authors have employed Chebyshev polynomials as a basis for $\Reg{}$ in spherical problems; see \cite{livermore_jones_worland_2007,boyd2011comparing,Hollerbach+2013,Marti-benchmark,townsend_etal_2016}. It is possible to use Chebyshev polynomials and achieve reasonably good results. A Chebyshev scheme has the advantage of a fast transform. We discuss this more in the next section. The difficulty for  Chebyshev methods is that it takes a many degrees of freedom just to represent the leading-order $r^{\ell}$ in $\Reg{}$ for $\ell\gg 1$. Depending on the error requirements, it may not require fully $\ell$ degrees of freedom, but it will not be easy. 

For example, it is possible to represent $r^{100}$ to within 1\% point-wise absolute error with Chebyshev polynomials of degree less than approximately 25. But the truncation will not allow for any of the analytic properties of the original function; \ie $r^{100} / r = r^{99}$, and $d r^{100}/dr = 100 r^{99}$. The Chebyshev truncation will give completely incorrect results under differentiation and division until the expansion degree reaches exactly 100. And even then, numerical imprecision can still lead to serious errors in the analytical properties. 

To solve the spherical Bessel problem, we first define $\vartheta \ \equiv \ r\, d/dr$. We  then need to multiply through by $r^{2}$ to avoid the singularity at $r=0$. That is, 
\Beq
\left[\, \vartheta ( \vartheta + 1 ) - \ell( \ell+ 1) + r^{2} \kappa^{2}\,\right] f_{\ell}(r) \ = \  0.   \label{regular SBFP}
\Eeq
It is advantageous to take the reflection symmetry into account in the Chebyshev expansion, but not to include an $r^{\ell}$ prefactor, which leads to severe ill conditioning \cite{boyd2011comparing}. Therefore expand,
\Beq
f_{\ell}(r) \ = \ r^{\gamma}\, \sum_{n=0}^{N}  T_{n}(2r^{2}-1)\, F_{n} , \quad \mathrm{where} \quad \gamma \ = \  \ell \! \!  \mod  2.   \label{Chebyshev expansions}
\Eeq
Under the action of the scaled derivative operator, $\vartheta$, 
\Beq
r \frac{d}{dr}  \left[ r^{\gamma}\, T_{n}(z) \right]  \ = \ r^{\gamma} \, \left[ 2 (1+z)T_{n}'(z) + \gamma \, T_{n}(z) \right] , \quad \mathrm{where} \quad z = 2r^{2} - 1. 
\Eeq
After posing the problem, sparse matrix construction requires using the following relationships between Chebyshev polynomials of the first and second kind,
\Beq
T_{n}'(z) \ = \ n U_{n-1}(z), \quad 2 T_{n}(z) \ = \ U_{n}(z) - U_{n-2}(z), \quad 2 z T_{n}(z) \ = \ T_{n+1}(z) + T_{n-1}(z) \label{Chebyshev relations}
\Eeq
\Eq{Chebyshev relations} are the same as the differentiation, conversion, multiplication operations we present in \S4, only between Jacobi polynomials with parameters $\pm 1/2$. \Eq{Chebyshev relations} also forms the foundation of all sparse Chebyshev methods that have come before. Most particularly the ``quasi-inverse'' method of Julien and Watson, and the ``ultra-spherical'' method of Olver and Townsend \cite{julien_watson_2009,olver2013fast}, both of which provided much inspiration for the methods in V16 and the current work. 

After using the properties of Chebyshev polynomials, the matrix version of \eq{regular SBFP} contains 9 matrix bands along the diagonal, and a dense boundary condition row. It is also possible to enforce the boundary conditions with basis recombination, and thereby increase the matrix bandwidth by two. The system is sparse, and straightforward to solve numerically. 

The biggest drawback of the Chebyshev method is the difficulty of working with vector and tensor systems. \Eq{regular SBFP} already takes more special intervention than in \S5.1. Everything  quickly becomes much more complicated  for higher-rank systems. It is also not possible to forward differentiate (\eg \textit{grad}), an expansion of the form in \eq{Chebyshev expansions}, which requires dividing by $r$. Townsend, Wilber, \& Wright \cite{townsend_etal_2016,wilber_etal_2017} took steps to mitigate this problem by filtering the first few monomials after solving. 

In many cases, a Chebyshev method is sufficient for reasonable results at modest resolution. But we can see from the error plots in Figure 1 that the results for more specialised methods can be orders of magnitude better.

\subsection{Grid representation and spectral transforms}

After finding the solution to a given system of equations in terms of spectral coefficients, it is natural to consider how to represent the solution on a spatial grid. Specifically, how can we compute the projections in \eq{projection}? The obvious answer is to use Gaussian quadrature. We implement this method for the calculation of nonlinear terms in Part-II. We use methods similar to those outlined in \cite{townsend_etal_2015}. Given practical and efficient transform methods for our bases, there are also some interesting of new alternatives on the horizon.

We discuss the issue of spectral transforms at length in V16. An unfortunate aspect of general Jacobi polynomials is the lack of a fast spectral transform. Only the Chebyshev polynomials have an unambiguous $\mathcal{O}(N \log N)$ transform. All other bases technically require $\mathcal{O}(N^{2})$ operations for the transforms. In V16, we demonstrated that number of spectral modes must get to a couple hundred before the asymptotically $\mathcal{O}(N \log N)$ beats matrix-based transform on modern SIMD hardware.  Boyd and Yu (2011) \cite{boyd2011comparing} also discuss the advantages and disadvantages of using a geometrically adapted method, versus a Chebyshev-based method with a fast transform. They also stress that the advantages of a fast transform are not as great as sometimes believed. For the 3D ball, the timing and scaling results for matrix transforms are identical to those in the disk. However, Jacobi-derived bases in three dimensions allow for additional advantages over 2D. 

Recall that,
\Beq
Q_{n}^{\alpha,\ell}(r) \ \propto \ P_{n}^{(\alpha, \ell+ 1/2)}(z)  
\Eeq
We know that we can convert each Jacobi index up one degree with a two-band matrix. This means that we can convert \textit{down} one degree by inverting a two-band matrix, which requires $\mathcal{O}(N)$ operations per conversion. If we pick $\alpha=-1/2$ as the primary basis, then we can convert our basis all the way to $P_{n}^{(-1/2, - 1/2)}(z)  \propto T_{n}(z)$ in $\mathcal{O}(N \ell)$ total operations. This is advantageous because we can use a fast transform on the Chebyshev basis. 

The downside to converting to a Chebyshev basis is that typically $\ell_{\max} \sim N$. There are big savings for small-$\ell$ values, but not so much for the high-degree modes. It is also not clear \textit{a priori} how numerically stable it is to convert down many degrees. But recent work by Slevinsky \cite{Slevinsky_SHT,slevinsky_2018,Slevinsky_etal} shows great promise in this regard for fast and stable algorithms with general Jacobi polynomial bases. Slevinsky showed how to use a divide-and-conquer butterfly method to compute the fast conversion from arbitrary Jacobi parameters to parameter values $-1 < a,b \le0$. For half-integer values of both parameters, it is then possible to use the Chebyshev transform after conversion. If the down-conversion does not land on $-1/2$, then Slevinsky \cite{slevinsky_2018} and others \cite{Hale_Townsend} recently demonstrated how to transform to the Chebyshev basis (from \eg Legendre) in a simple way based on asymptotics. While these methods are not fully exact, the accuracy can be made machine precision without compromising overall efficiency.  

Moreover, the results of Slevinsky and others allow for effective fast transforms for spin-weighted spherical harmonics. Although, fast spherical harmonics transforms are not a serious impediment to the creation of high-performance codes for multi-scale physical applications. Several highly tuned numerical transform libraries have proven quite effective in advanced applications \cite{SHTns,Libsharp}.

\section{Conclusions}

This paper derives a new algorithm for the accurate and efficient solution of partial differential equations involving scalars, vectors, and tensors in spherical polar coordinates.  We represent variables in a spectral basis derived from the geometric properties of the sphere.  The basis respects the tensorial nature of equations, and tensorial derivative operators (\ie \textsl{grad}, \textsl{div}, \textsl{curl}, etc.) are all represented as maximally sparse matrices. Furthermore, the basis automatically satisfies the regularity conditions at the coordinate singularities of spherical polar coordinates: the north and south poles, and the centre of the sphere.

We derive the basis in two steps.  First we consider the angular dependence (\S2).  Following the work of Phinney \& Buridge (1973) \cite{P&B73}, we use \textsl{spin-weighted} spherical harmonics.  These are a generalisation of the regular spherical harmonics (which correspond to a spin-weight of 0), which can also be used to represent vectors and tensors on a 2-sphere. To represent the angular variation of vectors and tensors, we first need to transform them into a basis of elements with definite spin-weight, via a unitary transformation. A new result for computations on the 2-sphere is the derivation of operators which represent the action of multiplication by $\cos(\theta)$ or $\sin(\theta)$, which is useful for the solution of equations whose coefficients have angular variation. Setting up a calculation for on the surface of a sphere is no more challenging than working with Fourier series. We have already deployed the above methods to study 2D active-matter turbulence on the surface of a sphere \cite{mickelin_etal_2018}. 

The second step derives a radial basis. The difficulty is that different components of tensors may have different behaviour near $r=0$. One key result of this paper is the derivation of an orthogonal transformation which takes a tensor represented with spin-weighted spherical harmonics, to a new basis in which each element is in a specific regularity space (\S3). This orthogonal transformation depends on the spherical harmonic degree $\ell$. A simple recursion relation allows the easy construction of the required transformations via \eq{Q-recursion}. This extends the work of James (1976) \cite{James76} which defined a different type of transformation based off Wigner 3-$j$ coefficients, and which was difficult to use for high-rank tensors.

A particular basis of weighted Jacobi polynomials represents each regularity space (\S4). Sparse tensorial derivative operators (\eg \textit{grad} and \textit{div}) map between different classes of weighted Jacobi polynomials. We derive a series of useful operators for the solution of partial differential equations, all based on the properties of Jacobi polynomials; which is another key result of this paper. These can be combined in a straightforward way to solve partial differential equations (\S5).

We also demonstrate how to construct different truncations via the $\tau$ method in \S5. We find that different residuals can lead to error properties that differ by many orders of magnitude. But we only demonstrate these differences for the simples relevant eigenvalue problem. In Part-II of this series, we show how our methods behave \textit{in situ}.  

Part-II describes the implementation of our algorithms as an extension of the Dedalus code framework\footnote{See: \url{http://dedalus-project.org}}. For demonstrating whole-code capability we solve several different partial differential equations relevant in astrophysical and geophysical fluid dynamics; written in tensorial form. For example, switching from solving a hydrodynamics problem to a magneto-hydrodynamics problem involves the addition of a few dozen of lines of code. We are able to demonstrate similar performance to other numerical methods which are specifically designed to magneto-hydrodynamics, despite our much more general approach. The algorithms derived in this paper thus have great promise for the solution of a wide range of problems involving partial differential equations in spherical polar coordinates.

As an overall theme, our work in this series of papers, and in V16, demonstrate how to construct very efficient and general numerical methods on domains with coordinate singularities. The approach of using the algebraic properties of Jacobi polynomials extends in a similar way to even more general domains.  Working out the details can get complicated. But the end result follows the same pattern and ends up with maximally sparse operators. With other collaborators, we have already shown how Jacobi polynomials have analogous properties on triangles \cite{triangle}. The corners of a triangle act the same as polar-coordinate singularities. It is necessary to link the Jacobi parameters to account for singular behaviour \cite{dunklxu,koornwinder}. The same situation holds for squares, and higher-dimensional polyhedra.

\bigskip

\textbf{ACKNOWLEDGMENTS}

GMV acknowledges support from the Australian Research Council, project number DE140101960.  DL is supported by a Hertz Foundation Fellowship, the National Science Foundation Graduate Research Fellowship under Grant No. DGE 1106400, a PCTS fellowship, and a Lyman Spitzer Jr. fellowship.

\bigskip

\textbf{Declarations of interest: none.}

\appendix

\section{Vector and tensor basis element notaion}

This appendix describes several aspects of tensor basis elements. Starting with the standard orthonormal coordinate basis vectors, we generate a \textit{spin-weighted} set of basis vectors,
\Beq
e_{0} \ \equiv \ e_{r},\quad 
\e{\pm} \ \equiv \ \frac{1}{\sqrt{2}}\left( \e{\theta} \mp i \e{\phi} \right). 
\Eeq
A unitary matrix transforms between the coordinate and spinor basis vectors,
\Beq
U \ \equiv \  \left[
\begin{array}{ccc}
 0 & \frac{1}{\sqrt{2}} &
   \frac{i}{\sqrt{2}} \\
 1 & 0 & 0 \\
 0 & \frac{1}{\sqrt{2}} &
   -\frac{i}{\sqrt{2}} \\
\end{array}
\right], \quad U^{-1} \ = \ U^{\dag}.
\Eeq
To transform higher-order tensors (of rank-$\mathfrak{r}$) between the coordinate and spinor bases, we take the Kronecker product of $U$ with itself $\mathfrak{r}$ times. 

The coordinate basis is orthonormal under the usual dot product
\Beq
 e_{r} \cdot e_{r}  \ = \  e_{\theta} \cdot e_{\theta}  \ = \ e_{\phi} \cdot e_{\phi} \ = \ 1,\quad  e_{r} \cdot e_{\theta}  \ = \  e_{r} \cdot e_{\phi}  \ = \ e_{\theta} \cdot e_{\phi} \ = \ 0
\Eeq
Therefore if $\sigma,\tau = -1, 0 , + 1$
\Beq
e_{\sigma} \cdot e_{\tau} =
\begin{cases} 1 & \mathrm{if} \ \sigma = - \tau \\ 
0 & \mathrm{otherwise} \end{cases}
\Eeq
This is the case for the standard \textit{dot product}. This is not a proper \textit{inner product} for the complex basis elements. We can take a couple of different approaches to defining dual (or contra-variant) vectors. We could simply use the complex conjugate, $e_{\sigma} \cdot e_{\tau}^{*} \ = \ e_{\sigma}^{*} \cdot e_{\tau} \ = \ \delta_{\sigma,\tau}$. We could alternatively define the contra-variant unit vector
\Beq
e^{\sigma} \ \equiv \ e_{-\sigma} \ = \ e_{\sigma}^{*}
\Eeq
Then one must take care of superscripts $e^{\sigma} \cdot e_{\tau} \ = \ \delta^{\sigma}_{\tau}$.  

The cross product for the coordinate basis implies,
\Beq
e_{0} \times e_{\pm} \ = \ \pm i\, e_{\pm}, \quad e_{+}\times e_{-} \ = \ i\, e_{0} 
\Eeq

We use multi-index notation for tensor components, 
\Beq
\sigma = \{\sigma_{1}, \ldots , \sigma_{\mathfrak{r}} \} 
\Eeq
We use a slightly different notation to denote a multi-index basis element, rather than an individual vector element. For tensor basis elements of a rank-$\mathfrak{r}$ tensor, we define 
\Beq
e(\sigma) \  \equiv \ e_{\sigma_{1}} \ldots e_{\sigma_{\mathfrak{r}}}  , \quad  \sigma_{i}\in \{-1,0,+1\}
\Eeq
The sum $\bar{\sigma}$ is called the spin-weight. 
\Beq
\bar{\sigma} \ \equiv \ \sum_{i=1}^{\mathfrak{r}} \sigma_{i}
\Eeq
This quantity, the spin-weight, arises repeatedly throughout the following analysis.
It is also useful to define the rank via
\Beq
|\sigma| \ \equiv \ \mathfrak{r}
\Eeq

A general rank-$\mathfrak{r}$ tensor therefore has the form
\Beq
\mathrm{T} \ = \ \sum_{\sigma} T^{\sigma} e(\sigma) , \quad \mathrm{where} \quad T^{\sigma}   \ =  \ T^{\sigma_{1}, \ldots , \sigma_{\mathfrak{r}}},
\Eeq
and each element $\sigma_{i} = -1, 0, +1$. 

For combining two sets of indices together,  
\Beq
\mu \sigma = \{\mu_{1}, \ldots,  \mu_{\mathfrak{r}}, \sigma_{1}, \ldots , \sigma_{\mathfrak{s}} \}.
\Eeq
Therefore
\Beq
e(\mu) e(\sigma) \ = \ e(\mu \sigma).
\Eeq
The multi-index notation incorporates scalars if $\sigma = \{ \}$, and vectors if $\sigma = {\sigma_{1}}$. In those cases
\Beq
e(\{\}) \ = \ 1, \quad e(\{\sigma\}) \ = \ e_{\sigma}.
\Eeq

We define the dual tensor index 
\Beq
\sigma^{*}  \ \equiv \ \left\{ - \sigma_{\mathfrak{r}}, \ldots, -\sigma_{1} \right\},
\label{sigma-star}
\Eeq
which is the negative and reverse of $\sigma$.
Under the action of vector-index contraction,
\Beq
e(\alpha^{*}) \cdot e(\beta) =  e_{-\alpha_{\mathfrak{r}}} \ldots e_{-\alpha_{1}} \cdot e_{\beta_{1}} \ldots e_{\beta_{\mathfrak{r}}} =  \delta(\alpha,\beta).
\Eeq

\section{Jacobi polynomial relationships}

\Beq
\label{Jacobi times-appendix}
\begin{array}{ccccc}
P_{n}^{(a,b)}(\z)             & = & \tfrac{n+a+b+1}{2n+a+b+1} P_{n}^{(a+1,b)}(\z)   & - & \tfrac{n+b}{2n+a+b+1} P_{n-1}^{(a+1,b)}(\z) \\ \\
\tfrac{1-\z}{2} P_{n}^{(a,b)}(\z) & = &   - \tfrac{n+1}{2n+a+b+1} P_{n+1}^{(a-1,b)}(\z) & + & \tfrac{n+a}{2n+a+b+1} P_{n}^{(a-1,b)}(\z)\\ \\
 P_{n}^{(a,b)}(\z) & = & \tfrac{n+a+b+1}{2n+a+b+1} P_{n}^{(a,b+1)}(\z)  & + & \tfrac{n+a}{2n+a+b+1} P_{n-1}^{(a,b+1)}(\z)  \\ \\ 
\tfrac{1+\z}{2} P_{n}^{(a,b)}(\z) &=&   \tfrac{n+1}{2n+a+b+1} P_{n+1}^{(a,b-1)}(\z) &+& \tfrac{n+b}{2n+a+b+1} P_{n}^{(a,b-1)}(\z).
\end{array}
\Eeq

\Beq
\frac{d}{d\z} P_{n}^{(a,b)}(\z) &=& \frac{n+a+b+1}{2} P_{n-1}^{(a+1,b+1)}(\z) \label{Jacobi Dp-appendix} \\
\left[b + (1+\z)\frac{d}{d\z}\right]P_{n}^{(a,b)}(\z) &=& (n+b) P_{n}^{(a+1,b-1)}(\z), \label{Jacobi Db-appendix} \\
\left[a - (1+\z)\frac{d}{d\z}\right]P_{n}^{(a,b)}(\z) &=& (n+a) P_{n}^{(a-1,b+1)}(\z), \label{Jacobi Da-appendix}
\Eeq

\Beq
P_{n}^{(a,b)}(\z=1) \ = \ {{n+a}\choose{a}} \ \equiv \ \frac{\Gamma(n+a+1)}{\Gamma(a+1)\, n!}.
\Eeq


\begin{thebibliography}{10}

\bibitem{Beyer_etal} 
{\sc F.~Beyer, B.~Daszuta, J.~Frauendiener, and B. Whale}, {\em Numerical evolutions of fields on the 2-sphere using a spectral method based on spin-weighted spherical harmonics}, Class. Quantum Grav. 31 (2014) 075019 

\bibitem{boyd_book} 
{\sc J.~P. Boyd}, {\em Chebyshev and Fourier Spectral Methods: Second Revised Edition}, Dover (2001)

\bibitem{boyd2011comparing} 
{\sc J.~P. Boyd and F.~Yu}, {\em Comparing seven spectral methods for
  interpolation and for solving the poisson equation in a disk: D+,D-
  polynomials, Logan--Shepp ridge polynomials, Chebyshev--Fourier series,
  cylindrical Robert functions, Bessel--Fourier expansions, square-to-disk
  conformal mapping and radial basis functions}, Journal of Computational
  Physics, 230 (2011), pp.~1408--1438.

\bibitem{coutsias_hagstrom_1996} 
{\sc E.~A.~Coutsias, T.~Hagstrom and D.~Torres}, {\em An efficient spectral method for ordinary
differential equations with rational function coefficients}, Mathematics of Computation, 65 (1996), pp. 611--635.

\bibitem{Charalambides_Waleffe_2008} 
{\sc M.Charalambides and F.~Waleffe}, {\em Gegenbauer tau methods with and without spurious eigenvalues},  SIAM J. Numerical. Analysis  47 (2008), pp. 48--68.

\bibitem{Dawkins_etal_1998} 
{\sc P.~T.~Dawkins, S.~R.~Dunbar, and R.W.~Douglass},  {\em The Origin and Nature of Spurious Eigenvalues in the Spectral Tau Method}, Journal of Computational Physics, 147 (1998), pp.~441--462.

\bibitem{Daou_Ortiz} 
M.~K.~EI-Daou, and E.~L.~Ortiz,
{The Tau Method as an Analytic Tool in the Discussion of Equivalence Results Across Numerical Methods} Computing, 60, (1998), pp.~365--376

\bibitem{doha2002efficient} 
{\sc E.~H. Doha and W.~M. Abd-Elhameed}, {\em Efficient spectral-Galerkin
  algorithms for direct solution of second-order equations using ultraspherical
  polynomials}, SIAM Journal on Scientific Computing, 24 (2002), pp.~548--571.

\bibitem{doha2006efficient} 
{\sc E.~H. Doha and A.~H.~Bhrawy}, {\em Efficient spectral-Galerkin algorithms for direct solution for second-order differential equations using Jacobi polynomials}, Journal of Numerical Algorithms, 42 (2006),  pp.~137--164

\bibitem{doha2009efficient} 
{\sc E.~Doha and W.~Abd-Elhameed}, {\em Efficient spectral
  ultraspherical-dual-Petrov--Galerkin algorithms for the direct solution of
  $(2n+ 1)$th-order linear differential equations}, Mathematics and Computers in
  Simulation, 79 (2009), pp.~3221--3242.

\bibitem{dunklxu} 
{\sc C.~F. Dunkl and Y.~Xu}, {\em Orthogonal polynomials of several variables,
  Second Edition}, Cambridge University Press, 2014.

\bibitem{Eastwood_Tod_1982} 
{\sc M.~Eastwood, and P.~Tod}, {\em Edth -- a differential operator on the sphere},
Math. Proc. Camb. Phil. Soc. 92 (1982),  pp.~317--330.

\bibitem{gardner_etal_1989} 
{\sc D.~R.~Gardner, S.~A. Trogden, and R.~W.~Douglass}, {\em A modified tau spectral method that eliminates spurious eigenvalues}, Journal of Computational Physics, 80 (1989), pp.~137--167.

\bibitem{G_S_56} 
{\sc I.~M.~Gelfand, and Z.~Y.~Shapiro},  {\em Representations of the group of rotations in three-dimensional space and their applications} , Am.~Math.~SOC~Transl., 2 (1956) pp.~207--316.

\bibitem{greengard_1991} 
{\sc L.~Greengard}
{\em Spectral integration and two-point boundary value problems} SIAM Journal of Numerical Analysis, 28 (1991), pp.~1071--1080

\bibitem{Hale_Townsend} 
{\sc N.~Hale. and A.~Townsend}, {\em A Fast, Simple, and Stable Chebyshev--Legendre Transform Using an Asymptotic Formula}, SIAM J. Sci. Comput., 36 (2014), pp.~A148--A167

\bibitem{Hollerbach+2013} 
{\sc R. Hollerbach,C.~Nore, P.~Marti,  S.~Vantieghem, F.~Luddens, and J.~L\'{e}orat}, {\em Parity-breaking flows in precessing spherical containers} Physical Review E 87 (2013), 053020

\bibitem{Homma}
{\sc Y.~Homma}, {\em Bochner--Weitzenb\"{o}ck formulas and curvature actions on Riemannian manifolds}, Trans. Amer. Math. Soc. 358 (2006),  87--114

\bibitem{James76} 
{\sc  R.~W.~James}, {\em New Tensor Spherical Harmonics, for Application to the Partial Differential Equations of Mathematical Physics} Phil.~Trans. Royal Society of London,  281 (1976), pp.~195--221

\bibitem{julien_watson_2009} 
{\sc K.~Julien, M.~Watson}, {\em Efficient multi-dimensional solution of PDEs using Chebyshev spectral methods}, Journal of Computational Physics, 228 (2009), pp.~1480--1503.

\bibitem{koornwinder} 
{\sc T.~Koornwinder}, {\em Two-variable analogues of the classical orthogonal
  polynomials}, Theory and applications of special functions,  (1975),
  pp.~435--495.

\bibitem{Kostelec} 
{\sc P.~J.~Kostelec, D.K.~Maslen, D.~M.~Healy, and D.~N.~Rockmore}
{\em Computational harmonic analysis for tensor fields on the two-sphere}, Journal of Computational Physics 162 (2000) pp.~514--535 

\bibitem{Lanczos38} 
{\sc C. Lanczos}, {\em Trigonometric interpolation of empirical and analytic functions}, J. Math. \& Physics , 17 (1938) pp.~123--199

\bibitem{Livermore} 
{\sc P.~W.~Livermore},  {\em Galerkin orthogonal polynomials} , Journal of Computational Physics, 229 (2010), pp~2046--2060

\bibitem{li2010optimal} 
{\sc K.~Li, P.~W. Livermore, and A.~Jackson}, {\em An optimal Galerkin scheme
  to solve the kinematic dynamo eigenvalue problem in a full sphere}, Journal
  of Computational Physics, 229 (2010), pp.~8666--8683.

\bibitem{Livermore_Ierley_2010} 
{\sc P.~W. Livermore, and G.~R.~Ierley}, {\em Quasi-$L^{p}$ norm orthogonal Galerkin expansions in sums of Jacobi polynomials}, Numer. Algor. 54 (2010) pp.~533--569

\bibitem{livermore_jones_worland_2007} 
{\sc P.~W. Livermore, C.~A.~Jones and S.~J.~Worland}, {\em Spectral radial basis functions for full sphere computations}, Journal of Computational Physics, 227 (2007), pp. 1209--1224

\bibitem{Marti_Jackson} 
{\sc P.~Marti, and A.~Jackson}, {\em A fully spectral methodology for magnetohydrodynamic calculations in a whole sphere} Journal of Computational Physics 305 (2016) pp.~403--422.

\bibitem{Marti-benchmark} 
{\sc P.~Marti, N.~Schaeffer, R.~Hollerbach, D.~C\'{e}bron, C.~Nore, F.~Luddens, J.-L.~Guermond, J.~Aubert, S.~Takehiro, Y.~Sasaki, Y.-Y.~Hayashi, R. Simitev, F.~Busse, S.~Vantieghem, and A.~Jackson},  {\em Full sphere hydrodynamic and dynamo benchmarks}, Geophys. J. Int. 197 (2014) , pp.~119--134.

\bibitem{matsushima1995spectral} 
{\sc T.~Matsushima and P.~Marcus}, {\em A spectral method for polar
  coordinates}, Journal of Computational Physics, 120 (1995), pp.~365--374.

\bibitem{mcfadden_etal_1990} 
{\sc G.~B. McFadden, B~T. Murray, and R~F.~Boisvert}, {\em Elimination of spurious eigenvalues in the Chebyshev tau spectral method}, Journal
  of Computational Physics, 91 (1990), pp.~228--239.

\bibitem{mickelin_etal_2018}
{O.~Mickelin, J~Stomka, K.~J.~Burns, D.~Lecoanet, G.M.~Vasil, L.~M.~Faria, and J.~Dunkel} {\em Anomalous chained turbulence in actively driven flows on spheres}, Phys. Rev. Lett., Accepted 21 March 2018

\bibitem{muite_2010} 
{\sc B.~K.~Muite}
{\em A numerical comparison of Chebyshev methods for solving fourth order semilinear initial boundary value problems}, Journal of Computational and Applied Mathematics, 234 (2010), pp. 317--342

\bibitem{N_P_66} 
{\sc E.~T.~Newman, and R.~Penrose}, {\em Note on the Bondi--Metzner--Sachs group}, J. Math. Phys. 7 (1966) pp.~863--870.

\bibitem{DLMF} 
{\sc F.~W.~J. Olver, D.~W. Lozier, R.~F. Boisvert, and C.~W. Clark}, {\em {NIST
  Handbook of Mathematical Functions}}, Cambridge University Press, 2010.

\bibitem{olver2013fast} 
{\sc S.~Olver and A.~Townsend}, {\em A fast and well-conditioned spectral
  method}, SIAM Review, 55 (2013), pp.~462--489.

\bibitem{triangle} 
{\sc S.~Olver, A.~Townsend, \& G.~M.~Vasil} {\em Recurrence relations for a family of orthogonal polynomials on a triangle}  \url{https://arxiv.org/pdf/1801.09099.pdf}

\bibitem{olver2014practical} 
{\sc S.~Olver and A.~Townsend}, {\em A practical framework for
  infinite-dimensional linear algebra}, in HPTCDL, 2014, pp.~57--62.

\bibitem{Ortiz69} 
{\sc E.~L.~Ortiz}, {\em The tau method}, SIAM J. Numer. Anal., 6 (1969) pp.~480--492

\bibitem{orszag_1974} 
{\sc S.A.~Orszag}, {\em Fourier series on spheres}, Monthly Weather Review, 102 (1974), pp.~56--75.

\bibitem{P&B73} 
{\sc R.~A.~Phinney, and R.~Buridge}, {\em Representation of the Elastic-Gravitational Excitation of a Spherical Earth Model by Generalized Spherical Harmonics}, Geophys J.~R.~asfr.~SOC. 34 (1973) , pp.~451--487.

\bibitem{Libsharp} 
{\sc M. Reinecke and D. S. Seljebotn},
{\em Libsharp --- spherical harmonic transforms revisited}
A\&A (2013) 554, A112.

\bibitem{Rubinstein}  
{\sc R.~Rubinstein, S.~Kurien, and C.~Cambon}, {\em Scalar and tensor spherical harmonics expansion of the velocity correlation in homogeneous anisotropic turbulence} Journal of Turbulence, 16 (2015), pp.~1058--1075

\bibitem{Slevinsky_SHT}  
{\sc R.~M.~Slevinsky}, {\em Fast and backward stable transforms between spherical harmonic expansions and bivariate Fourier series}, {Applied and Computational Harmonic Analysis} (2017)

\bibitem{slevinsky_2018}  
{\sc R.~M.~Slevinsky}, {\em On the use of Hahn's asymptotic formula and stabilized recurrence for a fast, simple, and stable Chebyshev--Jacobi transform}, 38 (2018) IMA Journal of Numerical Analysis, pp.~102--124

\bibitem{Slevinsky_etal} 
{\sc R.~M.~Slevinsky, H.~Montanelli, Q.~Du}, {\em A spectral method for nonlocal diffusion operators on the sphere} \url{https://arxiv.org/abs/1801.04902}

\bibitem{sakai2009application} 
{\sc T.~Sakai and L.~Redekopp}, {\em An application of one-sided Jacobi
  polynomials for spectral modelling of vector fields in polar coordinates},
  Journal of Computational Physics, 228 (2009), pp.~7069--7085.

\bibitem{SHTns} 
{\sc N.~Schaeffer}, {\em Efficient spherical harmonic transforms aimed at pseudospectral numerical simulations} Geochemistry Geophysics Geosystems 14 (2013),  pp.~751--758.

\bibitem{townsend_etal_2015} 
{\sc A. Townsend}, {\em The race for high order Gauss-Legendre quadrature}, in SIAM News, March 2015

\bibitem{townsend_etal_2016} 
{\sc A.~Townsend, H.~Wilber, \& G.~Wright}, {\em Computing with functions in spherical and polar geometries I. The sphere}. SISC, 38 (2016), C403-C425

\bibitem{trefethenspectralmethods} 
{\sc L.~N. Trefethen}, {\em Spectral methods in MATLAB}, Siam, 10, (2000.) 

\bibitem{Kronecker} 
{\sc C.~F.~ Van Loan}, {\em The ubiquitous Kronecker product}, J Computational \& Applied Mathematics, 123 (2000) pp.~85--100

\bibitem{V16} 
{\sc G.~M.~Vasil, K.~J. Burns, D.~Lecoanet, S.~Olver,
B.~P.~Brown, J.~S. Oishi}, {\em Tensor calculus in polar coordinates using Jacobi polynomials}, Journal of Computational Physics, 35 (2016), pp.~53--73.

\bibitem{viswanath_2015} 
{\sc D.~Viswanath}, {\em Spectral integration of linear boundary value problems}, Journal of Computational and Applied Mathematics, 290 (2015), pp.~159--173

\bibitem{wilber_etal_2017}
{\sc H.~Wilber, A.~Townsend  \& G.~B.~Wright}, {\em Computing with functions in spherical and polar geometries II. The disk}, SISC, 39 (2017), C238-C262

\bibitem{zernike} 
{\sc F.~Zernike}, {\em Beugungstheorie des schneidenverfahrens und seiner
  verbesserten form, der phasenkontrast-methode}, Physica, 1,  (1934),
  pp.~689--704.

\bibitem{zebib_1984} 
{\sc A.~Zebib}, {\em A Chebyshev method for the solution of boundary value problems}, Journal of Computational Physics,  53 (1984), pp.~443--455.

\end{thebibliography}
\end{document}